\documentclass[aos,MSNbibl,nameyear,dvips]{arximspdf}
\usepackage[section]{algorithm}
\usepackage{dcolumn}
\usepackage{algorithmic}
\usepackage{mathrsfs}

\doi{10.1214/12-AOS1039} 
\volume{40}
\issue{5}
\pubyear{2012}
\firstpage{2359}
\lastpage{2388}

\makeatletter

\newcolumntype{d}[1]{D{.}{.}{#1}}

\newproclaim{remarks}{Remarks}
\newproclaim{remark}{Remark}

\newproclaim{assumption}{Assumption}

\newproclaim{method}{Method}

\renewcommand{\widehat}{\hat}
\renewcommand{\widetilde}{\tilde}

\newcommand{\rright}{\right}
\newcommand{\lleft}{\left}

\newtheorem{theorem}{Theorem}
\newtheorem{lemma}[theorem]{Lemma}

\newcommand{\RR}{\mathbb R}
\newcommand{\EE}{\mathbb{E}}
\newcommand{\PP}{\mathbb{P}}

\makeatother

\begin{document}
\begin{frontmatter}

\title{Joint variable and rank selection for parsimonious estimation
of high-dimensional matrices}
\runtitle{Joint variable and rank reduction}

\begin{aug}
\author[A]{\fnms{Florentina} \snm{Bunea}\thanksref{t1}\ead[label=e1]{fb238@cornell.edu}},
\author[B]{\fnms{Yiyuan} \snm{She}\thanksref{t2}\ead[label=e2]{yshe@stat.fsu.edu}}
\and
\author[C]{\fnms{Marten H.} \snm{Wegkamp}\corref{}\thanksref{t1}\ead[label=e3]{marten.wegkamp@cornell.edu}}
\runauthor{F. Bunea, Y. She and M. H. Wegkamp}
\affiliation{Cornell University, Florida State University and Cornell
University}
\address[A]{F. Bunea\\
Department of Statistical Science\\
Cornell University\\
1198 Comstock Hall\\
Ithaca, New York 14853-3801\\
USA\\
\printead{e1}}
\address[B]{Y. She\\
Department of Statistics\\
Florida State University \\
117 N. Woodward Ave\\
Tallahassee, Florida 32306-4330\\
USA\\
\printead{e2}}
\address[C]{M. H. Wegkamp\\
Department of Mathematics\\
and\\
Department of Statistical Science\\
Cornell University\\
1194 Comstock Hall\\
Ithaca, New York 14853-3801\\
USA\\
\printead{e3}} 
\end{aug}

\thankstext{t1}{Supported in part by NSF Grant DMS-10-07444.}
\thankstext{t2}{Supported in part by NSF Grant CCF-1116447.}

\received{\smonth{10} \syear{2011}}
\revised{\smonth{4} \syear{2012}}

%
\begin{abstract}
We propose dimension reduction methods for sparse, high-dimen\-sional
multivariate response regression models. Both the number of responses
and that of the predictors may exceed the sample size. Sometimes viewed
as complementary, predictor selection and rank reduction are the most
popular strategies for obtaining lower-dimen\-sional approximations of
the parameter matrix in such models. We show in this article that
important gains in prediction accuracy can be obtained by considering
them jointly. We motivate a new class of sparse multivariate regression
models, in which the coefficient matrix has low rank \textit{and} zero
rows or can be well approximated by such a matrix. Next, we introduce
estimators that are based on penalized least squares, with novel
penalties that impose simultaneous row and rank restrictions on the
coefficient matrix. We prove that these estimators indeed adapt to the
unknown matrix sparsity and have fast rates of convergence. We support
our theoretical results with an extensive simulation study and two data
analyses.
\end{abstract}

%

\end{frontmatter}

\section{Introduction}\label{sec1}
\label{secintro}

The multivariate response regression model
%
\begin{equation}
\label{mod}Y = X A + E
\end{equation}
postulates a linear relationship between $Y$, the $m \times n$ matrix
containing measurements on $n$ responses for $m$ subjects, and $X$, the
$m \times p$ matrix of measurements on $p$ predictor variables, of rank
$q$. The term $E$ is an unobserved $m \times n$ matrix with independent
$N(0,\sigma^2)$ entries. The unknown $p \times n$ coefficient matrix
$A$ of unknown rank $r$ needs to be estimated. If we use (\ref{mod}) to
model complex data sets, with a high number of responses and
predictors, the number of unknowns can quickly exceed the sample size
$m$, but the situation need not be hopeless for the following reason.
Let $r$ denote the rank of $A$ and $J$ denote the index set of the
nonzero rows of $A$ and $|J|$ its cardinality. Counting the parameters
in the singular value decomposition of~$A$, we observe that in fact
only $r(n+|J|-r)$ free parameters need to be estimated, and this can be
substantially lower than the sample size~$m$. Furthermore, as we can
always reduce $X$ of rank $q$ to an $m\times q$ matrix with $q$
independent columns in $\RR^m$ that span the same space as the columns
of $X$, we can always assume that $|J| \leq q$. If $A$ is of full rank
with no zero rows, then the total number of parameters to be estimated
reverts back to $nq$. If either, or both, $q$ and $n$ are large, more
parsimonious models have to be proposed. Among the possible choices,
two are particularly popular.

The first class consists of
\textit{rank sparse} or \textit{rank deficient models}, which postulate
either that $A$ has low rank or that it can be well approximated by a
low rank matrix.
Methods tailored to rank sparsity seek adaptive rank $k$ approximations
of the coefficient matrix $A$. Then, one only needs to estimate $ k (q
+ n - k)$ parameters, which can be substantially less than $nq$ for low
values of $k$.

The second class of models reflects the belief that $|J|$ is smaller
than $q$, and we will call them \textit{row sparse models}. Methods
that adapt to \textit{row sparsity} belong to the variable selection
class, as explained in Section~\ref{sec1.1} below. The effective number of
parameters of such models is $|J|n$. This number is smaller than the
unrestricted~$nq$, but may be higher than $r(|J| + n - r)$, especially
if the rank of $A$ is low.

This discussion underlines the need for introducing and studying
another class of models, that embodies both sparsity constraints on $A$
simultaneously. In this work we introduce \textit{row and rank sparse}
models, and suggest and analyze new methods that combine the strengths
of the existing dimension reduction techniques. We propose penalized
least squares methods, with new penalties tailored to adaptive and
optimal estimation in the \textit{row and rank sparse} model (\ref{mod}).
The rest of the article is organized as follows.

We introduce in Section~\ref{sec2.1} a product-type penalty
that imposes simultaneously rank and row sparsity restrictions on the
coefficient matrix. It generalizes both AIC-type penalties developed
for variable selection in univariate response regression models as well
as the rank penalty of \citet{BuneaSheWegkamp} for low rank estimation
in multivariate response models. The purpose
of the resulting method is twofold. First, we prove in Theorem \ref
{thmeen} of Section~\ref{sec2.1} that the resulting estimators of $A$ adapt to
both types of sparsity, row and rank, under no conditions on the design
matrix. Their rates of convergence coincide with the existing minimax
rates in the literature, up to a logarithmic term; cf. \citet{kolt}.
Second, we show in Theorem~\ref{thmtwee} that this method can also be
employed for selecting among competing estimators from a large finite list.
This is of particular interest for selecting among estimates of
different ranks and sparsity patterns, possibly obtained via different
methods. The results of Section~\ref{sec2.1} hold for any values of $m, n$ and
$p$ and, in particular, both $n$ and $p$ can grow with $m$, but
computing the estimator analyzed in Theorem~\ref{thmeen} requires an
exhaustive search over the class of all possible models, the size of
which is exponential in $p$, and this becomes computationally
prohibitive if $p > 20$.

To address the computational issue, we propose two other methods in
Section~\ref{sec2.2}. The crucial ingredient
of both methods is the selection of predictors in multivariate response
regression models
\textit{under rank restrictions}. We define and analyze this core
procedure in Section~\ref{sec2.2}, and describe a computationally efficient
algorithm in Section~\ref{sec3.1}. By combining this method with two different
ways of selecting the rank adaptively we obtain two estimators of $A$.
Both are computable in high dimensions, and both achieve the rates
discussed in Section~\ref{sec2.1}, up to a $\log(p)$ factor, under different,
mild assumptions.
We also compare the theoretical advantages of these new methods over a
simple two-stage procedure in which one first selects the predictors
and then reduces the rank. We
illustrate the practical differences via a simulation study in Section
\ref{sec3.2}. We then use our methods for the analysis, presented in Section
\ref{secapps}, of two data sets arising in machine learning and cognitive
neuroscience, respectively. The proofs of our results are collected in
the \hyperref[app]{Appendix}.

\subsection{Background}\label{sec1.1}

Before we discuss our methods, we give an overview of existing
procedures of adaptive estimation in (\ref{mod}), that adapt to either
rank or row sparsity, but not both. We also present a comparison of
target rates under various sparsity assumptions on the coefficient
matrix $A$ in model (\ref{mod}).

Reduced rank estimation of $A$ in (\ref{mod}) and the immediate
extensions to principal components analysis (PCA) and canonical
correlation analysis (CCA) are perhaps the most popular ways of
achieving dimension reduction of multivariate data. They have become a
standard tool in time series [\citet{brill}], econometrics
[\citet{ReinVelu}] and machine learning [\citet{izenbook}],
to name just
a few areas. The literature on low rank regression estimation of $A$
dates back to \citet{And51}. The model is known as \textit{reduced-rank
regression} (RRR) [\citet{izenbook}] and, until recently, it had only
been studied theoretically from an asymptotic perspective, in a large
sample size regime. We refer to \citet{ReinVelu} for a historical
development and references, and to \citet{izenbook} for a large number
of applications and extensions. Very recently, a number of works
proposed penalized least squares estimators. For penalties proportional
to the nuclear norm, we refer to \citet{YuanNNP}, \citet{CandPlan},
\citet{NegahbanWainwright11}, \citet{Rohde11}. For penalties
proportional
to the rank, we refer to \citet{BuneaSheWegkamp} and \citet
{giraud}. Both
types of estimators are computationally efficient, even if $\max(n,p)>
m$, and both achieve, adaptively, the rate of convergence $(q+n)r$
which, under suitable regularity conditions, is the optimal minimax
rate in (\ref{mod}) under rank sparsity; see, for example,
\citet{Rohde11} for lower bound calculations.

To explain the other notion of sparsity, note that removing predictor
$X_j$ from model (\ref{mod}) is equivalent with setting the $j$th row
in $A$ to zero.
Since vectorizing both sides of model (\ref{mod}) yields a univariate
response regression model, we can view the rows of $A$ as groups of
coefficients in the transformed model. We can set them to zero by any
group selection method developed for univariate response regression
models in high dimensions such as the Group Lasso [\citet{Yuan}],
GLASSO for later reference.
The
optimal minimax rate in (\ref{mod}) under row sparsity is proportional
to $|J|n + |J| \log(p/|J|)$, again under suitable regularity
conditions; see
\citet{lounici-2010} and \citet{Weihuang}.

Despite these very recent advances, adaptive low rank estimation in
(\ref{mod}), based on a \textit{reduced} set of predictors, has not been
investigated either theoretically or practically.
For ease of reference, Table~\ref{tabl1} contains a rate comparison between
optimal prediction error rates achievable by variable selection
(GLASSO), low rank estimation (RSC and NNP) and our new joint rank and
row selection ({JRRS}) methods, respectively.

%
\begin{table}
\tablewidth=185pt
\caption{Oracle rate comparison between JRRS, GLASSO and RSC.
The sample size and dimension parameters $m, p, n, q, r, |J|$ satisfy
$q\leq m \wedge p $, $r \leq n \wedge|J|$, $|J| \leq q$}\label{tabl1}
\begin{tabular*}{\tablewidth}{@{\extracolsep{\fill}}r l@{}}
\hline
GLASSO: & $n |J|+ |J|\log(p)$ \\
RSC or NNP: & $nr+q r$ \\
JRRS: & $ nr + |J|r \log(p / |J|)$\\
\hline
\end{tabular*}
\end{table}

The table reveals that if $n\ge q$, the rates of the RSC, NNP and JRRS
are dominated by $nr$, regardless of $J$, while if $n<q$, the new class
of methods can provide substantial rate improvements over the existing
methods, especially when the rank $r$ is low.

\section{Adaptation to row and rank sparsity: Estimation procedures and
oracle inequalities}\label{sec2}

\subsection{The single-stage joint rank and row selection estimator}\label{sec2.1}
\label{onestep}

In this section we modify the rank selection criterion (RSC) introduced
in \citet{BuneaSheWegkamp}
to accommodate variable selection.
We propose our single-stage joint rank and row selection (JRRS) estimator
%
\begin{equation}
\label{step2} \widehat B = \mathop{\arg\min}_{B } \bigl\{\|Y - XB
\|^2_F + \operatorname{pen}(B) \bigr\},
\end{equation}
also denoted by JRRS1,
with penalty term
%
\begin{equation}\label{optpen-mult}
\operatorname{pen}(B) = c \sigma^2 r(B) \biggl\{ 2n + \log(2e)
\bigl|J(B)\bigr| +\bigl|J(B)\bigr| \log\biggl(\frac{ep}{|J(B)|} \biggr) \biggr\}.
\end{equation}
The penalty is essentially proportional to the number of parameters in
a model with fewer predictors $J(B)$ and of reduced rank $r(B)$. Here
$c > 3$ is a numerical constant,
$J(B)$ is the set of indices of nonzero rows, $r(B)$ is the rank of a
generic $p\times n$ matrix $B$ and the squared Frobenius norm of a
generic matrix $M$ is denoted by $\|M\|_F^2$ and is equal to the sum of
the squared entries of $M$.

If $\widehat B$ is computed by minimizing over \textit{all} $p \times n$
matrices $B$, then Theorem~\ref{thmeen} stated below shows that it
adapts optimally to the unknown row and rank sparsity of $A$: the mean
squared error of $\widehat B$ coincides with that of optimal estimators
of rank $r$ and with $|J|$ nonzero rows, had these values been known
prior to estimation. However, the construction of $\widehat B$ does not
utilize knowledge of either $r$ or $J$, hence the term adaptive. The
minimax lower bounds for this model can be obtained by an immediate modification
of Theorem 5 in \citet{kolt}. Our single-stage JRRS estimator
$\widehat
B$ given in (\ref{step2}) above achieves the lower bound, up to a log
factor, under no restrictions on the design $X$, rank $r$ or dimensions $m,n,p$.

\begin{theorem}\label{thmeen}
The single-stage JRRS estimator $\widehat B$ in (\ref{step2}) using
pen$(B)$ in (\ref{optpen-mult}) with $c=12$\setcounter{footnote}{2}\footnote{Our proof shows
that we may take $c>3$, at the cost of increasing numerical constants
in the right-hand side of the oracle inequality.}
satisfies
\[
\EE \bigl[ \| XA-X\widehat B\|_F^2 \bigr] \le 10 \|
XA-XB\|_F^2 + 8 \operatorname{pen}(B) + 768 n
\sigma^2 \exp(-n/2)
\]
for any $B$ with $r(B)\ge1$.
In particular, if $r(A)\ge1$,
\[
\EE \bigl[ \| XA-X\widehat B\|_F^2 \bigr] \lesssim
\sigma^2 r(A) \biggl\{ n+ \bigl|J(A)\bigr| \log\biggl( \frac{ p }{
|J(A)| } \biggr)
\biggr\}.
\]
Here and elsewhere $\lesssim$ means that the inequality holds up to
multiplicative numerical constants.
\end{theorem}

The proof of Theorem~\ref{thmeen} remains valid if the matrices we
select from depend on the data. Thus, our procedure
can be used for selecting from any countable list of random matrices
of different ranks and with different sparsity patterns.
We will make essential use of this fact in the next section.

\begin{theorem}\label{thmtwee}
For any collection of (random) nonzero matrices $B_1, B_2,\ldots, $
the single-stage JRRS estimator
%
\begin{equation}
\label{step3} \widetilde B = {\mathop{\arg\min}_{B_j }}\|Y -
XB_j\|^2_F + \operatorname{pen}(B_j)
\end{equation}
with $c=12$
satisfies
\begin{eqnarray*}
\EE \bigl[ \| XA-X\widetilde B\|^2_F \bigr] &\le&
\inf_{j} \bigl\{ 10 \EE\bigl[ \| XA-XB_j\|_F^2\bigr]
+ 8\EE\bigl[ \operatorname {pen}(B_j) \bigr] \bigr\}
\\
&&{}+ 768 n\sigma^2 \exp(-n/2).
\end{eqnarray*}
\end{theorem}

\subsection{Two-step joint rank and row selection estimators}\label{sec2.2}
The computational complexity of the single-stage JRRS estimator (\ref
{step2}) is owed to the component of the penalty
term proportional to $J(B)$, which is responsible for row selection.
The existence of this term in (\ref{optpen-mult}) forces complete
enumeration of the model space. We address this problem by proposing a
convex relaxation
$\| B\|_{2,1}$ of this component. Here $\|B\|_{2,1}=\sum_{j = 1}^{p}\|
b_j\|_2$ is the sum of the Euclidean norms $\| b_j\|_2 $ of the rows
$b_j$ of $B$.
In this section we propose two alternatives, each a two-step JRRS
procedure and each building on the following core estimator.

\subsubsection{Rank-constrained predictor selection}\label{sec2.2.1}
We define our rank-con-\break strained row-sparse estimators $\widehat B_{k}$
of $A$ as
%
\begin{equation}\label{grpsel-constrRank}
\widehat{B}_{k} = \mathop{\arg\min}_{r(B) \leq k} \bigl\{ \| Y
- XB\|_F^2 + 2\lambda\| B\|_{2,1} \bigr\}.
\end{equation}
Here $\lambda$ is a tuning parameter and the minimization is over all
$p\times n$ matrices $B$ of rank less than or equal to (a fixed) $k$.
A computationally efficient numerical algorithm for solving this
minimization problem is given in Section~\ref{secalg}. Clearly, for
$k=q$, there is no rank restriction in (\ref{grpsel-constrRank}) and
the resulting estimator is the GLASSO estimator; for $\lambda=0$, we
obtain the reduced-rank regression estimator. Thus, the procedure
yielding the estimators $\widehat{B}_{k}$ of rank $k$ acts as a
synthesis of the two dimension reduction strategies, having each of
them as limiting points. We will refer to $\widehat B_{k}$ as the rank
constrained group lasso (RCGL) estimators.

Since this estimator is central to our procedures, we analyze it first.
We need the following mild assumption on $\Sigma= X'X /m$.
{\renewcommand{\theassumption}{$\mathfrak{A}$}
\begin{assumption}\label{assumA}
We say $\Sigma\in{\mathbb R}^{p\times p}$ satisfies condition
$\mathfrak{A}(I, \delta_I)$ for an index set $I\subseteq\{1,\ldots,p\}$
and positive number $\delta_I$, iff
%
\begin{equation}
\label{condition1} \operatorname{tr}\bigl(M'\Sigma M\bigr) \ge
\delta_I \sum_{ i \in I} \|m_i
\|_2^2
\end{equation}
for all $p\times n$ matrices $M$ (with rows $m_i$) satisfying $ \sum_{
i \in I} \|m_i\|_2 \geq2\sum_{ i \in I^c} \|m_i\|_2$.
\end{assumption}}

\begin{remarks*}
(1) The constant 2 may be replaced by any constant larger than 1.

(2) Assumption~\ref{assumA} allows designs with $p>m$, and can be
seen as a version of the restricted eigenvalue condition in the
variable selection literature introduced in \citet{bickel-2008} and
analyzed in depth in \citet{Sara}.\vadjust{\goodbreak}

(3) A sufficient condition for (\ref{condition1}) is: there exists a
diagonal matrix $D$ with $D_{jj} = \delta_I$ for all $ j\in I$ and
$D_{jj}=0$ otherwise such that
$\Sigma- D$ is positive definite.\vspace*{-2pt}
\end{remarks*}

Let $ \lambda_1(\Sigma)$ denote the largest eigen-value of $\Sigma$ and
set the tuning parameter
%
\begin{equation}
\label{tune-lambda} \lambda= C \sigma\sqrt{\lambda_1(\Sigma) k m
\log(e p)}
\end{equation}
for some numerical constant $C > 0 $. Notice that $\lambda$ depends on
$k$, but we suppress this dependence in our notation.\vspace*{-2pt}
%
\begin{theorem}\label{thmdrie}
Let $\widehat B_{k}$ be the global minimizer of (\ref
{grpsel-constrRank}) corresponding to $\lambda$ in (\ref{tune-lambda})
with $C$ large enough. Then, we have
%
\begin{eqnarray}
\label{twee} \EE \bigl[ \| X \widehat B_{k} - X A
\|_F^2 \bigr]
&\lesssim&\| X B - X A \|_F^2 \nonumber\\[-8pt]\\[-8pt]
&&{}+ k \sigma^2
\biggl\{n + \biggl(1+ \frac
{\lambda_1(\Sigma) }{\delta_{J(B)}} \biggr)\bigl|J(B)\bigr|\log(p) \biggr\}
\nonumber
\end{eqnarray}
for any $ p\times n$ matrix $B$ with $1\le r(B)\leq k$, $|J(B)|$
nonzero rows, provided $\Sigma$ satisfies assumption $\mathfrak
{A}(J(B), \delta_{J(B)})$.\vspace*{-2pt}
\end{theorem}
\begin{remarks*}
(1)
The term $\{ n+ |J(B)|\log(p) \} k \sigma^2$ in (\ref{twee}) is
multiplied by a factor $\kappa= \lambda_1(\Sigma)/ \delta_{J(B)} $.
This factor can be viewed as a generalized condition number of the
matrix $\Sigma$.
If $\kappa$ stays bounded, Theorem~\ref{thmdrie}
shows that, within the class of row sparse matrices of fixed rank $k$,
the RCGL estimator
is row-sparsity adaptive, in that the best number of predictors does
not have to be specified prior to estimation.

(2)
It is interesting to contrast our estimator with the regular GLASSO
estimator $\widehat A$ that minimizes $\|Y-XB\|_F^2 + 2\lambda'
\| B\|_{2,1}$ over \textit{all}
$p\times n$ matrices~$B$.
Our choice (\ref{tune-lambda}) of the tuning parameter $\lambda$
markedly differs from the choice
proposed by \citet{lounici-2010} for the GLASSO estimator
$\widehat A$.
We need a different choice for $\lambda$ and a more refined analysis
since we minimize in (\ref{grpsel-constrRank}) over all $p\times n$
matrices $B$ of rank $r(B)\le k$.\vspace*{-2pt}
\end{remarks*}

\subsubsection{Adaptive rank-constrained predictor selection}\label{sec2.2.2}
We now develop theoretical properties of three methods, Method~\ref{method1}
(RSC$\to$RCGL), Method~\ref{method2} (RCGL$\to$JRRS1) and Method~\ref{method3}
(GLASSO$\to$RSC).
$JRRS1$ denotes the single-stage JRRS estimator of Section~\ref{sec2.1}.

Theorem~\ref{thmdrie} suggests that by complementing RCGL by a method
that estimates the rank consistently,
we could obtain row and rank optimal adaptive estimator. This is indeed
true.\vspace*{-2pt}

\begin{method}[(RSC$\to$RCGL)]\label{method1}
\begin{itemize}
\item
Use the rank\vspace*{1pt} selection criterion (RSC) of \citet{BuneaSheWegkamp} to
select $k=\widehat r$
as the number of singular values of $PY$ that exceed $ \sigma( \sqrt {2n}+\sqrt{2q})$. Here $P$ is the projection matrix on the space
spanned by $X$.\vadjust{\goodbreak}
\item Compute the rank constrained GLASSO estimator $\widehat B_{k}$
in (\ref{grpsel-constrRank}) above with $k=\widehat r$ to obtain the
final estimator $\widehat B^{(1)} = \widehat B_{\widehat r}$.\vspace*{-2pt}
\end{itemize}
\end{method}

This two-step estimator adapts to both rank and row sparsity, under two
additional, mild restrictions.
\begin{itemize}
\item[\textsc{Assumption} $\mathfrak{C1}$.]
$d_r(XA)> 2\sqrt{2}\sigma(\sqrt{n}+\sqrt{q})$.
\item[\textsc{Assumption} $\mathfrak{C2}$.]
$\log( \| XA\|_F ) \le(\sqrt{2}-1)^2 (n+q)/4$.
\end{itemize}

Assumption $\mathfrak{C1}$ only requires that the signal strength,
measured by $d_r(XA)$, the $r$th singular value of the matrix $XA$, be
larger than
the ``noise level'' $2\sigma\sqrt{2n+2q}$, otherwise its detection
would become problematic. The tightness of $\mathfrak{C1}$ is discussed
in detail in \citet{BuneaSheWegkamp}. Theorem 2 of that work
proves that
the correct rank will be selected with probability\footnote{Hence, if
$n+q$ is small compared to $m$, we suggest to replace $q$ by $q\log(m)$
in the threshold level in the definition of $\widehat r$, in $\mathfrak
{C1}$ and in $\mathfrak{C2}$; see the remark following Corollary 4 in
\citet{BuneaSheWegkamp}.}
$1-\exp\{-c_0(n+q)\}$ with $c_0=(\sqrt{2}-1)^2/2$.

Assumption $\mathfrak{C2}$ is technical and needed to guarantee that
the error due to selecting the rank is negligible compared to the rate
$nr +|J|r\log(p)$.\vspace*{-2pt}
%
\begin{theorem}\label{thmvijf}
Let $\Sigma$ satisfy $ \mathfrak{A}(J,\delta_J)$ with $J=J(A)\ne
\varnothing$, let $\lambda_1(\Sigma)/ \delta_{J}$ be bounded, and let
$\mathfrak{C1}$ and $\mathfrak{C2}$ hold.
Then the two-step JRRS estimator $\widehat B^{(1)}$ with
$\lambda$ set according to (\ref{tune-lambda}) with $C$ large enough
satisfies
\[
\EE \bigl[ \bigl\| X \widehat B^{(1)} - X A \bigr\|_F^2
\bigr] \lesssim nr + |J|r \log(p).\vspace*{-2pt}
\]
\end{theorem}

Hence, $\widehat B^{(1)} $ is \textit{row and rank} adaptive, and achieves
the same optimal rate, up to a $\log(p)$ factor, as the \textit{row and
rank} adaptive $\widehat B$ studied in Theorem~\ref{thmeen} above.
While Theorem~\ref{thmeen} is proved under no restrictions on the
design, we view the mild conditions
of Theorem~\ref{thmvijf} as a small price to pay for the computational
efficiency of
$\widehat B^{(1)} $ relative to that of $\widehat B$ in (\ref{step2}).
The practical choice of the threshold $2\sigma\sqrt{2n+2q}$ in the
initial step of our procedure can be done either by
replacing $\sigma^2$ by an estimator, as suggested and analyzed
theoretically in Section 2.4 of \citet{BuneaSheWegkamp}, or by
cross-validation. The latter is valid in this context for consistent
rank selection, as the minimum squared error of rank restricted
estimators in (\ref{mod}) is achieved for the true rank, as discussed
in detail in \citet{BuneaSheWegkamp}.

We now present an alternative adaptive method that is more
computationally involved than Method~\ref{method1}, as it involves a search
over a two-dimensional grid, but its analysis does not require
$\mathfrak{C1}$ and $\mathfrak{C2}$.\vspace*{-2pt}

\begin{method}[(RCGL$\to$JRRS1)]\label{method2}
\begin{itemize}
\item Pre-specify a grid $\Lambda$ of values for $\lambda$ and use
(\ref{grpsel-constrRank}) to construct the class $\mathcal{B}= \{
\widehat{B}_{k, \lambda}\dvtx  1\le k \le q,
\lambda\in\Lambda\}$.\vadjust{\goodbreak}
\item
Compute
$\widehat B^{(2)} = \mathop{\arg\min}_{B \in\mathcal{B}} \{\|Y
- XB\|^2_F + \operatorname{pen}(B) \}$,
with pen$(B)$ defined in (\ref{optpen-mult}) above.
\end{itemize}
\end{method}

We have the same conclusion as for Method~\ref{method1}:
%
\begin{theorem}\label{thmvier}
Provided
$\Sigma$ satisfies condition $\mathfrak{A}(J,\delta_J)$ with
$J=J(A)\ne
\varnothing$,
$\lambda_1(\Sigma)/\delta_{J} $ is bounded, and
$\Lambda$ contains $\lambda$ in (\ref{tune-lambda}) for some $C$ large
enough, we have
\[
\EE \bigl[ \bigl\| X\widehat B^{(2)}- XA\bigr\|^2_F
\bigr] \lesssim nr + |J| \log(p) r.
\]
\end{theorem}

We see that $\widehat B^{(2)}$ has the same rate as $\widehat B^{(1)}$,
under condition $\mathfrak{A}$ on the design only.
Our simulation studies in Section~\ref{sec4} indicate that the numerical results
of Methods~\ref{method1} and~\ref{method2} are comparable.
\begin{remark*}
A perhaps more canonical two-stage procedure is as follows:

\begin{method}[(GLASSO$\to$RSC)]\label{method3}
\begin{itemize}
\item Select the predictors via the GLASSO.
\item
Use the rank selection criterion (RSC) of \citet{BuneaSheWegkamp} to
construct an adaptive estimator, of reduced rank, based only on the
selected predictors.
\end{itemize}
It is clear that as soon as we have selected the predictors
consistently in the first step, selecting consistently the rank in the
second step and then
proving row and rank sparsity of the resulting estimator will follow
straightforward
from existing results, for instance, Theorem 7 in \citet
{BuneaSheWegkamp}. Although this is a natural path to follow, there is
an important caveat to consider: the sufficient conditions under which
this two-step process yields adaptive (to row and rank sparsity)
estimators include the conditions under which the GLASSO yields
consistent group selection. These conditions are
in the spirit of those given in \citet{FLORI}, for the Lasso, and
involve the mutual coherence condition on $\Sigma= (\Sigma_{ij})_{ 1
\leq i, j \leq p}$, which postulates that the off-diagonal elements of
$\Sigma$ be small. Specifically, for the GLASSO, the restriction
becomes $|\Sigma_{ij}| \leq1/(7\alpha|J|)$, for some $\alpha>1$ [cf.
\citet{lounici-2010}], if it is coupled with the condition that
$n^{-1/2} \|a_j\| \geq C_1 m^{-1/2} [1 + C_2n^{-1/2} \log p]^{1/2}$.
Here $\| a_j \|$ is the Euclidean norm of the $j$th row vector of $A$,
and $C_1$ and $C_2$ are constants.
For designs for which $\Sigma$ is even closer to the identity matrix,
in that
$|\Sigma_{ij}| \leq1/(14\alpha|J|n)$, the condition on the minimum
size of detectable coefficients can be relaxed to $n^{-1/2} \|a_j\|
\geq C_1 m^{-1/2} [1 + C_2n^{-1} \log p]^{1/2}$; see Corollary 5.2 in
\citet{lounici-2010}. Our Theorems~\ref{thmvijf} and~\ref{thmvier}
require substantially weaker assumptions on the design.\looseness=1
\end{method}
\end{remark*}

\section{Computational issues and numerical performance comparison}\label{sec3}

\subsection{A computational algorithm for the RCGL-estimator}\label{sec3.1}
\label{secalg}
In this section we design an algorithm for minimizing
%
\begin{equation}\label{grpsel-constrRank-prob}
F(B;\lambda):= \tfrac{1}{2} \| Y - XB\|_F^2 +
\lambda\| B\|_{2,1}\vadjust{\goodbreak}
\end{equation}
over all $p\times n$ matrices $B$ of rank less than or equal to $k$.
Recall that by solving this problem we provide a way of performing
\textit{rank-constrained variable selection} in model (\ref{mod}).
Directly solving the nonconvex constrained minimization problem for $B$
in (\ref{grpsel-constrRank-prob}) may be difficult. One way of
surmounting this difficulty is to write $B=SV'$, with $V$ being
orthogonal. Then the rank constrained group lasso (RCGL) optimization
problem is equivalent to finding
%
\begin{equation}\label{grpsel-constrRank2}
(\widehat S, \widehat V) = \mathop{\arg\min}_{S\in\mathbb
{R}^{p\times k}, V \in\mathbb{O}^{n\times k}}
\frac{1}{2} \bigl\| Y - X S V' \bigr\|_F^2 +
\lambda\| S \|_{2,1},
\end{equation}
where the minimum is taken over all orthogonal $n\times k$ matrices $V$
and all $p\times k$ matrices $S$.
With a slight abuse of notation, we still denote the objective function
in (\ref{grpsel-constrRank2}) by $F(S, V; \lambda)$.
We propose the following iterative optimization procedure.

\begin{algorithm}[t]
\renewcommand{\thealgorithm}{$\mathscr A$}
\begin{algorithmic}
\caption{$\quad$ Rank constrained group lasso (RCGL) computation
\strut
\label{algS-RRR}}
\STATE{\textbf{given}} $1\leq k \leq m \wedge p \wedge n$, $\lambda\geq
0$, $V_{k,\lambda}^{(0)} \in\mathbb{O}^{n\times k}$ (say the first $k$
columns of $I_{n\times n}$ or the first $k$ right singular vectors of $Y$)
\STATE$j\gets0$, $\mathrm{converged}\gets\mathrm{FALSE}$
\WHILE{not converged}
\STATE(a) $S_{k,\lambda}^{(j+1)} \gets\arg\min_{S \in\mathbb
{R}^{p\times k} } \frac{1}{2} \| Y V_{k,\lambda}^{(j)} - X S \|_F^2 +
\lambda\| S \|_{2,1}$.
\STATE(b) Let $W \gets Y' X S_{k,\lambda}^{(j+1)} \in\mathbb
{R}^{n\times k}$ and perform SVD: $W= U_w D_w V_w'$.
\STATE(c) $V_{k,\lambda}^{(j+1)} \gets U_w V_w'$
\STATE(d) $B_{k,\lambda}^{(j+1)} \gets S_{k,\lambda}^{(j+1)} (
V_{k,\lambda}^{(j+1)})'$
\STATE(e) $\mathrm{converged}\gets| F(B_{k,\lambda}^{(j+1)};
\lambda
)- F(B_{k,\lambda}^{(j)};\lambda)|<\varepsilon$
\STATE(f) $j\gets j+1$
\ENDWHILE
\STATE{\textbf{deliver}} $\widehat B_{k,\lambda}= B_{k,\lambda}^{(j+1)}$,
$\widehat S_{k,\lambda}= S_{k,\lambda}^{(j+1)}$, $\widehat
V_{k,\lambda
}= V_{k,\lambda}^{(j+1)}$.
\end{algorithmic}
\end{algorithm}

The following theorem presents a \textit{global} convergence analysis for
Algorithm~\ref{algS-RRR}, where \textit{global} in this context refers to
the fact that the algorithm converges for any initial
point.
%
\begin{theorem} \label{algstat}
Given $\lambda>0$ and an arbitrary starting point $V_{k,\lambda}^{(0)}
\in\mathbb{O}^{n\times k}$, let $(S_{k, \lambda}^{(j)}, V_{k,
\lambda
}^{(j)})$ ($j=1,2,\ldots$) be the sequence of iterates generated by
Algorithm~\ref{algS-RRR}.
The following two statements hold:
\begin{longlist}[(ii)]
\item[(i)] Any accumulation point of $(S_{k, \lambda}^{(j)}, V_{k,
\lambda}^{(j)})$ is a stationary point of $F$ and $F(S_{k, \lambda
}^{(j)}, V_{k, \lambda}^{(j)})$ converges monotonically to $F(S_{k,
\lambda}^{*}, V_{k, \lambda}^{*})$ for some stationary point $(S_{k,
\lambda}^{*}, V_{k, \lambda}^{*})$.
\item[(ii)] Suppose for any $(S, V)$ outside the local minimum set of $F$,
$ F(S, V)> \min_{\tilde S \in{\mathbb R}^{p\times k}} F(\tilde S, V)$.
Then, any accumulation point of $(S_{k, \lambda}^{(j)}, V_{k, \lambda
}^{(j)})$ is a local minimum of $F$ and $F(S_{k, \lambda}^{(j)}, V_{k,
\lambda}^{(j)})$ converges monotonically to $F(S_{k, \lambda}^{*},
V_{k, \lambda}^{*})$ for some local minimizer $(S_{k, \lambda}^{*},
V_{k, \lambda}^{*})$.
\end{longlist}
\end{theorem}
\begin{remarks*}
(1)
We run the algorithm to obtain a solution path, for each $(k, \lambda)$
in a two-dimensional grid or for a grid of $\lambda$ with $k$
determined by RSC. From the solution path, we get a series of candidate
estimates. Then the single stage JRRS (\ref{step3}) or other tuning
criteria can be used to select the optimal estimate.

(2) Our results are of the same type as those established for the
convergence of the EM algorithm [\citet{jeffEM}]. Algorithm
\ref{algS-RRR} can be viewed as a block coordinate descent method, but
the conclusion in Theorem~\ref{algstat} is stronger in some sense: the
guaranteed convergence to a stationary point (to be defined in Appendix
\ref {subsecglobalconv}) does \textit{not} require the uniqueness of
$S_{k,\lambda}^{(j+1)}$ in step (a) which is a crucial assumption in
the literature [see, e.g., \citet{Bert} and \citet{Tseng01}].

(3) Step (a) needs to solve a GLASSO optimization problem. To see this,
denoting the standard vectorization operator by $\mathrm{vec}$ and the
Kronecker product by $\otimes$, we rewrite $\| YV_k - X S_k\|_F^2/2 +
\lambda\| S_k\|_{2,1}$ as $ \|{\operatorname{vec}}( (Y V_k)') - ( X
\otimes I_k) \operatorname{vec}( S_k')\|_F^2/2 + \lambda\|S_k
\|_{2,1}$. Although this subproblem is convex, finding its global
minimum point can still be expensive for large data. Instead, one may
perform some low-cost thresholding for a few steps. Concretely, let $K$
be a constant satisfying $K > \|X\|_2^2/2$. Given $V \in{\mathbb
O}^{p\times k}$, define $T_V\dvtx  {\mathbb R}^{p\times k} \rightarrow
{\mathbb R}^{p\times k}$ as
%
\begin{equation}\label{grpsel-th}
T_V \circ S = \vec\Theta \biggl(\frac{1}{K}X' Y
V+ \biggl(I - \frac{1}{K} X' X\biggr) S;\frac{\lambda}{K}
\biggr)\qquad \forall S \in{\mathbb R}^{p\times k},
\end{equation}
where $\vec\Theta$ is a multivariate version of the soft-thresholding
operator $\Theta$. For any vector $a\in\mathbb R^{k}$, $\vec\Theta
(a;\lambda):=
a \Theta(\|a\|_2;\lambda)/\|a\|_2$ for $a\neq0$ and $0$ otherwise;
for any matrix $A\in{\mathbb R}^{p\times k}$ with $A= [ a_1\enskip
\cdots\enskip a_p  ]'$, $\vec\Theta(A;\lambda):= [ \vec\Theta
(a_1;\lambda) \enskip\cdots\enskip\vec\Theta(a_p;\lambda) ]'$.

We now replace step (a) in Algorithm~\ref{algS-RRR} by $S_{k,\lambda}^{(j+1)}\gets
T_{V_{k,\lambda}^{(j)}} \circ\cdots\circ T_{V_{k,\lambda}^{(j)}}
\circ S_{k,\lambda}^{(j)}$,
where the number of $T_{V_{k,\lambda}^{(j)}}$, denoted by $\alpha_j$,
satisfies $1\leq\alpha_j \leq M_{\mathrm{iter}}$ for some $M_{\mathrm{iter}}<\infty$
specified based on available computational resources. $\alpha_j$ need
not be equal.
This algorithm, denoted by $\mathscr A'$, offers more flexibility and
is more convenient than Algorithm~\ref{algS-RRR} in implementation.
Although at each iteration $S_{k,\lambda}^{(j+1)}$ is not uniquely
determined, a stronger global convergence result holds for $\mathscr
A'$.
\end{remarks*}
%
\begin{theorem} \label{algstat2}
Given $\lambda>0$ and an arbitrary starting point $V_{k,\lambda}^{(0)}
\in\mathbb{O}^{n\times k}$, let $(S_{k, \lambda}^{(j)}, V_{k,
\lambda
}^{(j)})$ ($j=1,2,\ldots$) be the sequence of iterates generated by
$\mathscr A'$. Then, any accumulation point of $(S_{k, \lambda}^{(j)},
V_{k, \lambda}^{(j)})$ is a coordinatewise minimum point (and a
stationary point) of $F$ and $F(S_{k, \lambda}^{(j)}, V_{k, \lambda
}^{(j)})$ converges monotonically to $F(S_{k, \lambda}^{*}, V_{k,
\lambda}^{*})$ for some coordinatewise minimum point $(S_{k, \lambda
}^{*}, V_{k, \lambda}^{*})$.
\end{theorem}

\subsection{Simulation studies}\label{sec3.2}
\label{secsimu}
The setup of our simulations is as follows:
\begin{itemize}
\item
The design matrix $X$ has i.i.d. rows $X_i$ from a multivariate normal
distribution $ \operatorname{MVN}(\bolds{0}, \Sigma)$, with
$\Sigma_{jk}=\rho^{|j-k|}$, $\rho> 0$, $1\leq j,k\leq p$.
\item
The coefficient matrix $A$ has the form
\[
A = \lleft[ %
\matrix{ A_1
\cr
O } %
\rright]=
\lleft[ %
\matrix{ b B_0 B_1
\cr
O }
\rright]
\]
with $b > 0$, $B_0$ a $J \times r$ matrix and $B_1$ a $r\times n$
matrix. All entries in $B_0$ and $B_1$ are i.i.d. $N(0, 1)$.
\item
The noise matrix $E$ has independent $N(0,1)$ entries. Let $E_i$ denote
its ith row.
\item
Each row $Y_i$ in $Y$ is then generated as $Y_i = X_{i}^{\prime}A+
E_i$, $1 \leq i \leq m$.
\end{itemize}
This setup contains many noisy features, but the relevant features lie
in a low-dimensional subspace. This structure
resembles many real world data sets; see our examples in Section \ref
{secapps}, where the low rank structure is inherent and, thus,
rank-constrained variable selection is desired.

We report two settings:
\begin{longlist}
\item[\textit{$p>m$ setup}:]
$m = 30$, $|J| = 15$, $p = 100$, $n = 10$, $r = 2$, $\rho=0.1$,
$\sigma^2=1$, $b=0.5,1$.
\item[\textit{$m>p$ setup}:]
$m = 100$, $|J| = 15$, $p = 25$, $n = 25$, $r = 5$, $\rho=0.1$,
$\sigma^2=1$,
$b=0.2, 0.4$.
\end{longlist}
Although we performed experiments in many other settings, say, with
$\rho=0.5$, we do not report all results, as the conclusions are
similar. The current setups show that variable selection, without
taking the rank information into consideration, may be suboptimal even
if the correlations between predictors are low.

We tested five methods: RSC, GLASSO, Method~\ref{method1}
(RSC$\to$RCGL), Meth\-od~\ref{method2} (RCGL$\to$JRRS1), and Method
\ref{method3} (GLASSO$\to$RSC), as described in Section~\ref{sec2.2}.
To minimize the influence of various parameter tuning strategies on our
performance comparison, we generated a large validation data set
(10,000 observations) to tune the parameter of each algorithm (with the
exception of Method~\ref{method2}) and we also generated another
independent data set of the same size as the test data to evaluate the
test error. Similar to the {LARS-OLS hybrid} [\citet{Efron2004}],
for each GLASSO and RCGL estimate, we computed the least squares
estimate restricted to the selected dimensions. We found that the
resulting (bias corrected) solution paths are more suitable for
parameter tuning. For\vspace*{1pt} Method~\ref{method2}, after getting
the (bias corrected) solution path, we set $c=3$ and $\sigma^2=1$ in
(\ref{optpen-mult}) to select the optimal $\widehat B^{(2)}$; in
contrast to the other two methods, no validation data is used for
tuning.\looseness=-1

Each model was simulated 50 times, and Tables~\ref{tabl2} and
\ref{tabl3} summarize our findings.
We evaluated the prediction accuracy of each estimator $\widehat A$ by
the mean squared error (MSE) $\| X A - X \widehat A\|_F^2/(mn)$ using
the test data at each run. Since the MSE histograms turned out to be
highly asymmetric, we computed
the $40\%$ trimmed-mean of MSEs as the goodness of fit of the
obtained model. This trimmed mean is more robust than the mean and more
stable than the median,
and it therefore allows for a more fair comparison between
methods.

\begin{table}
\tablewidth=235pt
\caption{Performance comparisons between GLASSO, RSC,
Methods \protect\ref{method1}, \protect\ref{method2} and \protect\ref{method3}
in the {$p>m$} experiment with $b=0.5, 1$, $|J|=15$ and $r=2$}
\label{tabl2}
\begin{tabular*}{\tablewidth}{@{\extracolsep{\fill}}l r r r r r@{}}
\hline
& \multicolumn{1}{c}{\textbf{MSE}} & \multicolumn{1}{c}{$\bolds{|\widehat J|}$}
& \multicolumn{1}{c}{$\bolds{\widehat R}$} & \multicolumn{1}{c}{\textbf{M}}
& \multicolumn{1}{c@{}}{\textbf{FA}} \\
\hline
\multicolumn{6}{@{}c@{}}{$b=0.5$}\\
[3pt]
\textit{GLASSO} & 206 & 10 & 10 & 53\% & 4\% \\
\textit{RSC} & 485 & 100 & 2 & 0\% & 100\% \\
[3pt]
\textit{Method}~\ref{method1} & 138 & 19 & 2 & 36\% & 10\% \\
\textit{Method}~\ref{method2} & 169 & 21 & 2 & 45\% & 7\% \\
\textit{Method}~\ref{method3} & 169 & 10 & 2 & 53\% & 4\% \\
[6pt]
\multicolumn{6}{@{}c@{}}{$b=1$}\\
[3pt]
\textit{GLASSO} & 511 & 14 & 10 & 41\% & 7\% \\
\textit{RSC} & 1905 & 100 & 2 & 0\% & 100\% \\
[3pt]
\textit{Method}~\ref{method1} & 363 & 21 & 2 & 31\% & 12\% \\
\textit{Method}~\ref{method2} & 434 & 25 & 2 & 30\% & 13\% \\
\textit{Method}~\ref{method3} & 402 & 14 & 2 & 41\% & 7\% \\
\hline
\end{tabular*}
\end{table}

\begin{table}
\tablewidth=235pt
\caption{Performance comparisons between GLASSO, RSC, Methods
\protect\ref{method1}, \protect\ref{method2} and \protect\ref{method3}
in the {$m>p$} experiment with $b=0.2, 0.4$, $|J|=15$ and $r=5$}
\label{tabl3}
\begin{tabular*}{\tablewidth}{@{\extracolsep{\fill}}l r r r r r@{}}
\hline
& \multicolumn{1}{c}{\textbf{MSE}} & \multicolumn{1}{c}{$\bolds{|\widehat J|}$}
& \multicolumn{1}{c}{$\bolds{\widehat R}$} & \multicolumn{1}{c}{\textbf{M}}
& \multicolumn{1}{c@{}}{\textbf{FA}} \\
\hline
\multicolumn{6}{c@{}}{$b=0.2$}\\
[3pt]
\textit{GLASSO} & 18.1 & 14 & 14 & 4\% & 1\% \\
\textit{RSC} & 11.9 & 25 & 5 & 0\% & 100\% \\
[3pt]
\textit{Method}~\ref{method1} & 8.3 & 15 & 5 & 0\% & 1\% \\
\textit{Method}~\ref{method2} & 8.3 & 15 & 5 & 0\% & 4\% \\
\textit{Method}~\ref{method3} & 8.9 & 14 & 5 & 4\% & 1\% \\
[6pt]
\multicolumn{6}{c@{}}{$b=0.4$}\\
[3pt]
\textit{GLASSO} & 17.7 & 15 & 15 & 0\% & 0\% \\
\textit{RSC} & 11.5 & 25 & 5 & 0\% & 100\% \\
[3pt]
\textit{Method}~\ref{method1} & 8.1 & 15 & 5 & 0\% & 0\% \\
\textit{Method}~\ref{method2} & 8.0 & 15 & 5 & 0\% & 1\% \\
\textit{Method}~\ref{method3} & 8.1 & 15 & 5 & 0\% & 0\% \\
\hline
\end{tabular*}
\end{table}

We also report the median number of predictors (denoted by $|\widehat
J|$) and median {rank} estimate (denoted by $\widehat R$) over all runs.
Estimators with small MSE and low $|\widehat J|$ and $\widehat R$ are
preferred from the point of view of statistical modeling.\looseness=-1

Finally, we provide the rates of nonincluded true variables (denoted by
$\mathrm{M}$ for misses) and the rates of incorrectly included variables
($\mathrm{FA}$ for false alarms). Ideally, both rates are low, especially
the M-rates, since we do not wish to discard relevant features.

We can draw the following conclusions from Tables~\ref{tabl2} and
\ref{tabl3}:
\begin{itemize}
\item
We see that straightforward variable selection via GLASSO often
severely misses some true features in the {$p>m$} setup as seen from
its high M numbers. RSC achieved good rank recovery, as expected, but,
by the definition of this estimator, it uses all $p$ variables.
Clearly both GLASSO and RSC alone are inferior to the three JRRS-type
methods (Methods~\ref{method1},~\ref{method2} and~\ref{method3}).

\item
Method~\ref{method1} dominates all other methods.
Its MSE results are impressive and confirm the rate improvement
established in Section~\ref{sec2.2} over the GLASSO and RSC. While
Method~\ref{method1} may not give exactly $|\widehat J| = |J|=15$, its M numbers indicate
that we did not miss many true features.

\item
Method~\ref{method2}, unlike Methods~\ref{method1} and~\ref{method3}, did not use
the large validation data for ideal parameter tuning, which explains
its slight inferiority relative to the other two methods. However, we
see that even without validation-based tuning, which may at times be
infeasible in practice, this method is a serious contender. It supports
the theoretical findings of Theorem~\ref{thmtwee} on the usage of the
penalty (\ref{optpen-mult}) for model comparison and tuning parameter
selection.

\item The performance of Method~\ref{method3} is inferior to that of
Method~\ref{method1} in the $p>m$ experiment, and comparable with both
Methods~\ref{method1} and~\ref{method2} in the $m>p$ experiment, when the
three methods have
essentially the same behavior.
\end{itemize}

In conclusion, we found that Method~\ref{method1} is the clear winner in
terms of performance as well as computational speed, among the
two-stage JRRS procedures we considered, and is particularly appealing
in the $m < p$ regime. In particular, it shows the advantage of the
novel penalty type which enforces simultaneous (row) sparsity and rank
reduction on the coefficient matrix. Method~\ref{method2} using penalty (\ref
{optpen-mult}) provides evidence of success of Theorem \ref
{thmtwee}.

\section{Applications}\label{sec4}\label{secapps}
In this section we apply Method~\ref{method1}, with its tuning parameters
chosen via cross-validation, to two real data sets from machine
learning and cognitive neuroscience.

\subsection*{Norwegian paper quality} These data were obtained from a
controlled experiment that was carried out at a paper factory in Norway
(Norske Skog, the world's second-largest producer of publication paper)
to uncover the effect of three control variables $X_1, X_2, X_3$ on the
quality of the paper which was measured by 13 response variables. Each
of the control variables $X_i$ takes values in $\{-1, 0, 1\}$. To
account for possible interactions and nonlinear effects, second order
terms were added to the set of predictors, yielding $X_1, X_2, X_3,
X_1^2, X_2^2, X_3^2, X_1 \cdot X_2, X_1\cdot X_3, X_2 \cdot X_3$ and
the intercept term. There were 29 observations with no missing values
made on all response and predictor variables. The Box--Behnken design of
the experiment and the resulting data are described in \citet{Aldrin96}
and \citet{izenbook}. Since neither the group penalty nor the rank
constraint is imposed on the intercept term, we always center the
responses and standardize the predictors in the training data (and
transform the validation/test data accordingly).

The data set can be downloaded from the website of \citet
{izenbook} and
its structure
clearly indicates that dimension reduction is possible, making it a
typical application for reduced
rank regression methods. The RSC method with
adaptive tuning, as described in \citet{BuneaSheWegkamp}, selected the
rank $\widehat r=3$.
This finding is consistent with \citet{Aldrin96}, who assessed the
performance of the rank 3 estimator by leave-one-out cross-validation
(LOOCV) and obtained a minimum LOOCV error (total squared error,
unscaled) of 326.2. We then employed the newly developed Method~\ref{method1}
to automatically determine the useful predictors and pursue the optimal
projections. Not surprisingly, the selected rank is still 3,
yielding 3 new scores, which are now constructed from only 6
of the original 9 predictors, with $X_1^2$, $X_1\cdot X_2$ and $X_2
\cdot X_3$
%
\begin{table}
\caption{Paper quality control: joint new feature construction from
the selected predictors with $\widehat r =3, |\widehat J|=6$. The
extracted components are ordered by their associated eigenvalues}
\label{tabl4}
\begin{tabular*}{\tablewidth}{@{\extracolsep{\fill}} l d{3.1} c d{2.4} d{2.4} d{2.4}
d{2.4} d{2.4}@{}}
\hline
\textbf{New scores} &
\multicolumn{1}{c}{\textbf{Eigenvalues}} &
\multicolumn{1}{c}{$\bolds{X_1}$} & \multicolumn{1}{c}{$\bolds{X_2}$} &
\multicolumn{1}{c}{$\bolds{X_3}$} & \multicolumn{1}{c}{$\bolds{X_2^2}$}
& \multicolumn{1}{c}{$\bolds{X_3^2}$} &
\multicolumn{1}{c@{}}{$\bolds{X_1 \cdot X_3}$}\\
\hline
1 & 112.4 & 1.9244 & -0.5288 & -1.2321 & -0.4443 & 1.3109 & 1.1898\\
2 & 40.7 & 0.8231 & -0.5937 & 0.9324 & 0.7819 & -0.1599 & -0.8536\\
3 & 24.9 & 0.2871 & 1.0336 & 0.4215 & 0.8365 & 0.4245 & 0.4677\\
\hline
\end{tabular*}
\end{table}
discarded, and only the variables from Table~\ref{tabl4}
selected. The tuning result was the same for 10-fold CV and LOOCV. The
minimum LOOCV error is now 304.5. We found no interaction effect
between $X_2$ and $(X_1, X_3)$, an interesting complement to Aldrin's
analysis. Table~\ref{tabl4} shows the construction weights
of the 3 new orthogonal score variables from the rank-3 RSC on the
selected set of variables.
They are ordered by an importance measure given by the associated
eigenvalues of $Y' X (X' X)^{-1} X' Y=:W$ [see \citet{ReinVelu} and
\citet{izenbook} for the explanation]. For instance, the first important
score variable (accounting for 57.5\% of the trace of~$W$) can be
roughly read as $2 X_1 -0.5 X_2 - X_3 - 0.5 X_2^2 + X_3^2 + X_1 X_3$,
or simply $2 X_1 + X_1 X_3 + 1_{X_3=-1}-1_{X_2=1}$. This can be used as
a concise summary predictor for all 13 response variables
simultaneously and it quantifies the effect of the design variables on
paper quality control.

\subsection*{Cognitive neuroimaging}

We present an analysis of the data set described in \citet
{tutorial} and
collected to investigate the effect of the HIV-infection on human
cognitive abilities. Neuro-cognitive performance is typically measured
via {correlated} neuro-cognitive indices (NCIs). This study employed $n
= 13$ NCIs, falling into five domains of attention/working memory,
speed of information processing, psychomotor abilities, executive
function, and learning and memory. These indices were measured for 62
HIV$+$ patients in the study. The set of explanatory variables was large
and contained: (a) clinical and demographic predictors and (b) brain
volumetric and diffusion tensor imaging (DTI) derived measures of
several white-matter regions of interest, such as fractional
anisotropy, mean diffusivity, axial diffusivity and radial diffusivity,
along with all volumetrics${}\times{}$DTI interactions. We refer to
\citet{tutorial} for details. The final model has $p=235$
predictors, much
greater than the sample size $m = 62$.
An initial analysis of this data set was performed using the
RSC to select a model of rank 1 and construct the corresponding new
predictive score. Although this is a massive reduction of the dimension
of the predictor space, all 235 initial predictors were involved in the
construction of the new score.

This leaves unanswered the important question as to what variables
(especially which DTI derived measures) are most predictive of
the neuro-cognitive changes in HIV$+$ patients. After standardizing the
predictors, we run Method~\ref{method1}.

We selected a model of rank 1 and constructed one new predictive score
but, very importantly, this score is a linear combination of \textit{only}
10 predictors that were selected from the original pool of 235.

When we set aside $30\%$ of the data as a
separate test set, and used the remaining to fit the model and tune the
regularization parameters, the mean squared error (MSE) of the RSC
estimate was $192.9$, while the MSE of the newly proposed method was
only $138.4$. Moreover, our analysis not only demonstrates the
existence of a strong association between the variable
\texttt{Education} and the neuro-cognitive abilities of HIV$+$ patients,
which had already been established by other means in the literature,
but also suggests, as a perhaps new finding, that the variable
fractional anisotropy at corpus callosum (\texttt{fa\_cc1}) stands out
among the very many DTI-derived measures, in terms of predictive
power.

\begin{appendix}\label{app}
\section*{Appendix}
For a generic index set $I\subseteq\{1,\ldots,p\}$, we define the
$m\times p$ matrix
$X_I$ as follows: its $i$th column coincides with that of $X$ if $i\in
I$, otherwise we set the entire column to zero.
Furthermore,
we
define $P_I$ as the projection matrix on the column space of $X_I$.

Since we favor transparent proofs, we did not attempt to optimize
various numerical constants.

\subsection{\texorpdfstring{Proof of Theorem \protect\ref{thmeen}}{Proof of Theorem 1}}\label{secA.1}
By the definition of $\widehat B$, for any $p\times n$ matrix $B$ with
$r(B)\ge1$, the inequality
\[
\| Y-X\widehat B\|_F^2 + \operatorname{pen}(\widehat B)
\le \| Y-XB\|_F^2+ \operatorname{pen}(B)
\]
holds.
This is equivalent with
%
\begin{eqnarray}
\label{begin} \| XA-X\widehat B\|_F^2 & \le& \| XA-XB
\|_F^2 + 2 \operatorname {pen}(B)
\nonumber\\[-8pt]\\[-8pt]
&&{} + \bigl\{ 2\langle E,X\widehat B-XB\rangle- \operatorname {pen}(\widehat B)-
\operatorname{pen}(B) \bigr\}.
\nonumber
\end{eqnarray}
We consider two complementary cases: $r(\widehat B) \ge1$ and
$r(\widehat B)=0$.

\textit{Case} 1: $r(\widehat B)\ge1$.
We write $J=J(B)$ and $\widehat J=J(\widehat B)$, and we note that
$J\cup\widehat J = J_1\cup J_2\cup J_3$ with $J_1 = \hat J^c \cap J$,
$J_2 = J \cap\hat J$, $J_3 = J^c \cap\hat J$.
Hence,
\[
\bigl\langle E,X(\widehat B- B)\bigr\rangle = \sum
_{k=1}^3 \bigl\langle E, X_{J_k} (\widehat
B-B) \bigr\rangle.
\]
The penalty term in (\ref{optpen-mult}) can be written as
\[
\operatorname{pen}(B) = c\sigma^2 r(B) \bigl\{ 2n + f\bigl(\bigl|J(B)\bigr|
\bigr) \bigr\}
\]
for the function
$f(x)= x \log(2e) + x \log(e p/x)$.
Since $f''(x)=-1/x$, $f(x)$ is concave for $x>0$, and we have
\[
f(x+y) \ge \frac12 f(2x) + \frac12 f(2y) = x+y+x\log\biggl(\frac
{ep}{x}
\biggr) + y\log\biggl(\frac{ep}{y}\biggr)
\]
for all $x,y>0$.
Consequently, writing $r_1=r(B)$, $ r_3=r(\widehat B)$ and $r_2=r_1+r_3$,
\begin{eqnarray*}
\operatorname{pen}(B) + \operatorname{pen}(\widehat B) &=& c
\sigma^2r_1 \bigl\{ 2n + f\bigl(|J_1|+|J_2|\bigr)
\bigr\} + c\sigma^2 r_3 \bigl\{ 2n + f\bigl(|J_2|+|J_3|\bigr)
\bigr\}
\\
&\geq& c\sigma^2\sum_{k=1}^3
r_k \biggl\{n+ |J_k| + |J_k| \log\biggl(
\frac{e
p}{|J_k|}\biggr) \biggr\}.
\end{eqnarray*}
This implies that
\begin{eqnarray*}
&& \bigl\{ 2\langle E,X\widehat B-XB\rangle- \operatorname {pen}(\widehat B)-
\operatorname{pen}(B) \bigr\}
\\
&&\qquad\le\sum_{k=1}^3 \biggl[ 2\bigl\langle
E, X_{J_k} (\widehat B-B) \bigr\rangle- c\sigma^2r_k
\biggl\{ n+|J_k|+|J_k| \log \biggl(\frac{e p}{|J_k|}
\biggr) \biggr\} \biggr].
\end{eqnarray*}
We define
\[
R = \max_{I} \biggl( d_1^2(P_{I}
E)- \frac{c}{2} \sigma^2 \biggl\{ n+| I |+|I| \log\biggl(
\frac{e p}{|I|}\biggr)\biggr\} \biggr)_{+}.
\]
In this proof,\footnote{A careful inspection reveals that we may take
any $c>3$, at the cost of larger constants elsewhere.} we set $c=12$.
Using the inequality
\begin{eqnarray*}
\bigl\langle E, X_{J_k} (\widehat B-B) \bigr\rangle& =& \bigl\langle
P_{J_k} E, X_{J_k} (\widehat B-B) \bigr\rangle
\\
& \le& d_1(P_{J_k}E) \sqrt{r_k} \bigl\|
X_{J_k} (\widehat B-B)\bigr\|_F
\end{eqnarray*}
and
the inequalities $2xy\le x^2/a + y^2 a$ and $(x+y)^2\le x^2 (1+a)+ y^2
(1+a)/a$ for all $x,y\in\RR$ with $a=2$, we further bound
\begin{eqnarray*}
&& \bigl\{ 2\langle E,X\widehat B-XB\rangle- \operatorname {pen}(\widehat B)-
\operatorname{pen}(B) \bigr\}
\\
&&\qquad\le\sum_{k=1}^3 \biggl[ \frac12 \bigl\|
X_{J_k}(\widehat B-B)\bigr\|_F^2 + 2 r_k
d_1^2(P_{J_k} E)\\
&&\quad\qquad\hspace*{18.2pt}{} - c\sigma^2
r_k \biggl\{ n+|J_k|+|J_k| \log \biggl(
\frac{e
p}{|J_k|}\biggr) \biggr\} \biggr]
\\
&&\qquad\le \frac12 \| X\widehat B-XB\|_F^2 + 6 \Bigl(
\max_k r_k\Bigr) R
\\
&&\qquad\le \frac34 \| X\widehat B- XA\|_F^2 + \frac32 \| XB
- XA\|_F^2 + 12n R,
\end{eqnarray*}
using $\max(r_1,r_2,r_3)=r_2\le2n$.
Now, (\ref{begin}) and the display above
yield
\begin{eqnarray*}
\tfrac14 \EE \bigl[ \|X\widehat B - XA\|_F^2 \bigr] &
\le& \tfrac52 \| XB-XA\|_F^2 + 2 \operatorname{pen}(B) +
12 n \EE[ R ]
\\
&\le&\tfrac52 \| XB-XA\|_F^2 + 2 \operatorname{pen}(B) +
192n\sigma^2 \exp(-n/2)
\end{eqnarray*}
using Lemma~\ref{maxbound} in the last inequality.
This concludes the first case.\vadjust{\goodbreak}

\textit{Case} 2: $r(\widehat B)=0$. Using the same reasoning as above, we
can argue that
\begin{eqnarray*}
\| XA \|_F^2 &=& \| XA - X\widehat B
\|^2_F
\\
& \le& \| XB -XA\|_F^2 + 2\operatorname{pen}(B) - 2\langle E, XB
\rangle- \operatorname{pen}(B)
\\
&\le&\tfrac52\| XB -XA\|_F^2 + 2\operatorname{pen}(B) + \tfrac34
\| XA\|_F^2+ R'
\end{eqnarray*}
with
\begin{eqnarray*}
R' &= & 2r(B) d_1^2(P_{J(B)}E) -
\operatorname{pen}(B)
\\
&\le& 2 r(B) \biggl\{ d_1^2(P_{J(B)} E) -
\frac{c}{2} \sigma^2\bigl( n + \bigl|J(B)\bigr| \bigr) \biggr\}.
\end{eqnarray*}
By Lemma 3 in \citet{BuneaSheWegkamp},
we have
$\EE[ d_1^2(P_{J(B)} E) ] \le2\sigma^2 (n+|J(B)|)$, so that
$\EE[ R'] \le0 $ for our choice $c=12$ above. Hence,
\[
\tfrac{1}{4} \EE \bigl[ \| XA - X\widehat B\|^2_F
\bigr] \le \tfrac 52\| XB -XA\|_F^2 + 2\operatorname{pen}(B),
\]
which concludes the second case, and our proof.
%
\begin{lemma} \label{maxbound}
We have
\[
\EE \biggl[ \max_{J} \biggl( d_1^2(
P_JE) - 6\sigma^2 \biggl\{ n+ |J|+ |J| \log \biggl(
\frac{e p}{|J| } \biggr) \biggr\} \biggr)_{+} \biggr] \le 16
\sigma^2 \exp(- n/2).
\]
\end{lemma}
\begin{pf}
Notice that $d_1^2(P_J E)$ has
the same distribution as $d_1^2(Z_k)$ for a $m\times k$ matrix $Z_k$ of
independent $N(0,\sigma^2)$ entries with $k=|J|$.
Consequently, for any $\nu_k\ge0$,
\[
\EE \Bigl[ \max_{J} \bigl( d_1^2(
P_JE) - \sigma^2 \nu_{|J|}^2
\bigr)_{+} \Bigr] \le \sum_{k=1}^q
\pmatrix{p
\cr
k} \EE \bigl[ \bigl( d_1^2(Z_k)
- \sigma^2 \nu_k^2 \bigr)_{+}
\bigr].
\]
We write $V_k= d_1(Z_k) /\sigma$ and $\mu_k=\sqrt{n}+\sqrt{k}$,
$1\le
k\le p$.
From Lemma 3 in \citet{BuneaSheWegkamp}, we have
\[
\PP\{V_k - \mu_k \ge t\} \le\exp\bigl(-t^2/2
\bigr)
\]
for all $t>0$.
Hence, for
$ \nu_k = \mu_k+\delta_k $ and $\delta_k > 0$, we obtain
\begin{eqnarray*}
\EE \bigl[ V_k^2 - \nu_k^2
\bigr]_{+} &\le& \EE \bigl[ V_k^21 \bigl\{
V_k^2 \ge\nu_k^2 \bigr\} \bigr]
\\
&\le& \int_{\mu_k+\delta_k}^\infty2x \PP\{ V_k \ge x\}
\,dx
\\
&\le& \int_{\delta_k}^\infty2(\mu_k+s) \exp
\bigl(-s^2/2\bigr) \,ds
\\
&=& 2 \biggl( 1+\frac{\mu_k}{\delta_k} \biggr) \exp\bigl(-{\delta_k}^2/2
\bigr)
\end{eqnarray*}
for $1\le k\le p$.
Choosing $\delta_k^2 = 2k \log(ep/k) + n+ k$, we find
\begin{eqnarray*}
\EE \Bigl[ \max_{J} \bigl( d_1^2(
P_JE) - \sigma^2 \nu_{|J|}^2
\bigr)_{+} \Bigr] &\le& \sigma^2 \sum
_{k=1}^q \frac{(e p)^k}{k^k} 2 \biggl(1+
\frac{\mu_k}{{\delta}_k} \biggr) \exp\bigl(-{\delta}_k^{2}/2
\bigr)
\\
&\le& 6 \sigma^2 \sum_{k=1}^q
\frac{(e p)^k}{k^k} \exp\bigl(-{\delta }_k^{2}/2\bigr)
\\
&=& 6 \sigma^2 \sum_{k=1}^q
\exp(-k/2)\exp(-n/2)
\\
& \le&16\sigma^2 \exp(-n/2).
\end{eqnarray*}
Since
$ \nu_k^2 \le2\delta_k^2 + 2\mu_k^2 \le4k\log(ep/k) + 6k+6n$, the
claim follows.
\end{pf}

\subsection{\texorpdfstring{Proof of Theorem \protect\ref{thmtwee}}{Proof of Theorem 2}}\label{secA.2}
By the same reasoning as in the proof of Theorem~\ref{thmeen}, we obtain
\begin{eqnarray*}
\frac14 \bigl[ \|X\widetilde B - XA\|_F^2 \bigr] &\le&
\frac52 \| XB_j-XA\|_F^2 + 2
\operatorname{pen}(B_j)
\\
&&{} + 12 n \max_{I} \biggl( d_1^2(P_{I}
E)- \frac{c}{2} \sigma^2 \biggl(n+|I|+ |I| \log
\frac{e p}{|I|}\biggr) \biggr)_{+}
\end{eqnarray*}
for any $B_j$, random or not, with $r(B_j)>0$. Taking expectations on
both sides of the inequality and applying Lemma~\ref{maxbound}
gives the result.

\subsection{\texorpdfstring{Proof of Theorem \protect\ref{thmdrie}}{Proof of Theorem 3}}\label{secA.3}
We denote the row vectors of the matrix $\widehat B_k$ by $\widehat
b_1,\ldots, \widehat b_p$, and we write
$\widehat J_k= \{ i\dvtx  \widehat b_i\ne0\}$ for the index set of nonzero rows.
Let $B$ be any matrix with row vectors $b_1,\ldots, b_p\in\RR^n$ with
$r(B)\le k$ and $J= \{ j\dvtx  b_j\ne0\}$ such that $\Sigma=X'X/m$
satisfies condition $\mathfrak{A}(J,\delta_J)$.
We use $\tilde X$ to denote $X_{\tilde J}$ with $\tilde{J} = \widehat
J_k \cup{ J}$ and
we write $\Delta^2=\| XB- X A\|_F^2$ and $\widehat\Delta^2_k= \|
X\widehat B_k - XA\|_F^2$.
In this notation, we have by the definition of $\widehat B_k$
\[
\widehat\Delta^2_k + 2\lambda\|\widehat B_{k}
\|_{2,1} \leq \Delta^2 +2\bigl\langle E, X (\widehat
B_k - B)\bigr\rangle+2\lambda\|B\|_{2,1}.
\]
This implies that
%
\begin{eqnarray}
\label{eqeen} && \widehat\Delta^2_k + 2\lambda\bigl\| (B-
\widehat B_k)_{J^c} \bigr\|_{2,1}
\nonumber\\[-8pt]\\[-8pt]
&&\qquad\le\Delta^2 + 2\bigl\langle E, \tilde X (\widehat
B_k-B)\bigr\rangle+ 2\lambda\bigl\| (B-\widehat B_k)_J
\bigr\|_{2,1}.
\nonumber
\end{eqnarray}
For
the second term on the right of (\ref{eqeen}), we note that
\begin{eqnarray*}
2\bigl\langle E, \tilde X (\widehat B_k-B)\bigr\rangle&=& 2\bigl\langle
P_{\tilde J} E, \tilde X (\widehat B_k-B)\bigr\rangle
\\
&\le& 2d_1(P_{\tilde J} E) \sqrt{2k} \bigl\| \tilde X(\widehat
B_k -B)\bigr\|_F
\\
&\le& 2d_1(P_{\tilde J} E) \sqrt{2k} (\widehat
\Delta_k + \Delta)
\\
&\le& 16 k d_1^2(P_{\tilde J} E) +
\tfrac{1}{4} \bigl( \widehat\Delta^2 + \Delta^2
\bigr)
\end{eqnarray*}
using $r(\widehat B_k -B)\le r(\widehat B_k )+ r(B) \le2k$ and the
inequality $2xy\le4 x^2+y^2/4 $ for all $x,y\in\RR$.
This bound and inequality (\ref{eqeen}) give the inequality
%
\begin{eqnarray}
\label{eqtwee} &&\tfrac34 \widehat\Delta_k^2 + 2\lambda\bigl\|
( \widehat B_k - B)_{J^c} \bigr\|_{2,1}
\nonumber\\[-8pt]\\[-8pt]
&&\qquad\le\tfrac54 \Delta^2 + 16 k d_1^2(P_{\widetilde J}
E) + 2 \lambda \bigl\| ( \widehat B_k - B)_{J}
\bigr\|_{2,1}.
\nonumber
\end{eqnarray}
For the remainder of the proof, we consider two complementary cases.

\textit{Case} 1. \textit{Assume that}
%
\begin{equation}
\label{case1} \tfrac54 \Delta^2 + 16 k d_1^2(P_{\widetilde J}
E) \le2 \lambda \bigl\| ( \widehat B_k - B)_{J}
\bigr\|_{2,1}.
\end{equation}
In this case, (\ref{eqtwee}) implies that
\[
\bigl\|(\widehat B_k - B)_{J^c} \bigr\|_{2,1} \le 2 \bigl\|
(\widehat B_k - B)_{J} \bigr\|_{2,1}.
\]
Since
$\Sigma$ satisfies $\mathfrak{A}(J, \delta_J)$, we have
\[
\| X\widehat B_k - X B\|_F^2 \ge m
\delta_J \sum_{j\in J} \| \widehat
b_i -b_i \|_2^2.
\]
This inequality and the inequality $2xy\le x^2/a + y^2 a$ applied to
$a= 8/\delta_J$ give
%
\begin{eqnarray}
\label{crap} 2\lambda\bigl\| (B-\widehat B_k)_J
\bigr\|_{2,1} &=& 2\lambda\sum_{i\in J} \|
b_i -\widehat b_i\|_2
\nonumber
\\
&\le& \frac{a}{m}\lambda^2 |J| + \frac{m}{ a } \sum
_{i\in J} \| b_i -\widehat b_i
\|_2^2
\nonumber
\\
&\le& \frac{a }{m}\lambda^2 |J| + \frac{1}{a \delta_J} \| X
\widehat B_k - XB\|_F^2
\\
&\le& \frac{ a }{m}\lambda^2 |J| + \frac{2}{a \delta_J} \bigl(
\widehat \Delta^2_k + \Delta^2\bigr)
\nonumber
\\
&=& \frac{8 \lambda^2 |J| }{m\delta_J} + \frac{1}{4} \bigl(\widehat \Delta_k^2
+\Delta^2\bigr).
\nonumber
\end{eqnarray}
After we combine (\ref{eqtwee}), (\ref{case1}) and (\ref{crap}), we obtain
\[
\frac34 \widehat\Delta_k^2 \le 4 \lambda \bigl\| ( \widehat
B_k - B)_{J} \bigr\|_{2,1} \le\frac{16 \lambda^2 |J|
}{m\delta_J} +
\frac{1}{2} \bigl(\widehat\Delta_k^2 +
\Delta^2\bigr),
\]
hence,
%
\begin{equation}
\label{UNO} \widehat\Delta_k^2 \le 2
\Delta^2 + \frac{64 \lambda^2 |J| }{m\delta_J} = 2 \Delta^2 + 64C k
\sigma^2 |J| \log(ep) \frac{ \lambda_1(\Sigma)}{\delta_J}
\end{equation}
for the choice $\lambda^2 = C \lambda_1(\Sigma) m k \sigma^2 \log(ep)$.
This concludes the first case.\vadjust{\goodbreak}

\textit{Case} 2. \textit{Assume that}
%
\begin{equation}
\label{case2} \tfrac54 \Delta^2 + 16 k d_1^2(P_{\widetilde J}
E) > 2 \lambda \bigl\| ( \widehat B_k - B)_{J}
\bigr\|_{2,1}.
\end{equation}
In this case, (\ref{eqtwee}) now gives
%
\begin{equation}
\label{drie} 3 \widehat\Delta_k^2 \le 10
\Delta^2 + 64 k d_1^2(P_{\widetilde J}E).
\end{equation}
By Lemma~\ref{maxbound}, we have
%
\begin{eqnarray}
\label{vier}\quad \EE\bigl[ d_1^2 (P_{\tilde J} E)
\bigr] &\le& 6\sigma^2 n + 16\sigma^2 \exp (-n/2) + 6
\sigma^2 \EE \biggl[ | \tilde J| + |\tilde J| \log\biggl(
\frac
{ep}{|\tilde J|}\biggr) \biggr]
\nonumber\\[-8pt]\\[-8pt]
&\le& 6\sigma^2 n + 16\sigma^2 \exp(-n/2) + 12
\sigma^2 \log(ep ) \bigl( |J| + \EE \bigl[ |\widehat J_k| \bigr]
\bigr).
\nonumber
\end{eqnarray}
We now bound $\EE[ |\widehat J_k|] $.
Since $r(\widehat B_k) \leq k$, we may write $\widehat B_k=\widehat S
\widehat V_k'$ for some $\widehat V_k$ with $k$ columns satisfying
$\widehat V_k' \widehat V_k =I_{k\times k}$.
Following the lines of argument in the proof of Lemma~\ref{algconv},
with $V$ fixed at $\widehat V_k$, $\widehat S$ is the (globally)
optimal solution to the convex problem of (\ref{subprob1}).
Let $x_i$ and $\widehat s_i$ for $1\le i\le p$ be the column vectors of
$X$ and $\widehat S'$, respectively.
Using the Karush--Kuhn--Tucker condition of $\widehat S' $, we obtain
%
\begin{equation}
\|\widehat b_i\| \neq0 \quad\Rightarrow\quad\|\widehat s_i\|
\neq0 \quad\Rightarrow\quad \bigl\| x_i' (X \widehat S - Y \widehat
V_k) \bigr\|_2 = \lambda,
\end{equation}
so that
\begin{eqnarray*}
|\widehat J_k | \lambda^2 &=& \sum
_{j\in\widehat J_k } \bigl\| x_j'(X\widehat S-X A
\widehat V_k - E\widehat V_k)\bigr\|_2^2
\\
&=& \sum_{j\in\widehat J_k } \bigl\| x_j'(X
\widehat S-X A \widehat V_k - P_{\widehat J_k } E\widehat
V_k)\bigr\|_2^2
\\
&\le& 2\sum_{j\in\widehat J_k } \bigl\| x_j'(X
\widehat S-X A \widehat V_k) \bigr\|_2^2 + 2 \sum
_{j\in\widehat J_k } \bigl\| x_j'
P_{\widehat J_k } E\widehat V_k\bigr\|_2^2
\\
&\le& 2\lambda_1\bigl(X'X\bigr) \|X\widehat
B_k - X A \|_F^2 + 2\lambda_1
\bigl(X'X\bigr) \| P_{\widehat J_k } E\widehat V_k
\|_F^2
\\
&\le& 2m \lambda_1(\Sigma) \bigl\{\widehat\Delta^2_k
+k d_1^2(P_{\widehat J_k }E) \bigr\}
\end{eqnarray*}
using $r(\widehat V_k)\le k$.
Taking expectations on both sides,
%
\begin{eqnarray}
\label{vijf} \EE\bigl[ |\widehat J_k | \bigr] &\le& \frac{2m \lambda_1(\Sigma) }{\lambda^2} \bigl\{
\EE\bigl[ \widehat\Delta^2_k\bigr] +k\EE\bigl[
d_1^2(P_{\widehat J }E)\bigr] \bigr\}
\nonumber\\[-8pt]\\[-8pt]
&\le& \frac{2m \lambda_1(\Sigma) }{\lambda^2} \bigl\{ \EE\bigl[ \widehat \Delta^2_k
\bigr] +k\EE\bigl[ d_1^2(P_{\tilde J }E)\bigr] \bigr
\}
\nonumber
\end{eqnarray}
since $\widehat J_k \subseteq\tilde J$.
After we combine (\ref{vier}) and (\ref{vijf}), we get
\begin{eqnarray*}
\EE \bigl[ d_1^2( P_{\tilde J} E) \bigr] &\le& 6
\sigma^2 n + 16\sigma^2 \exp(-n/2) + 12
\sigma^2 \log(ep ) |J|
\\
&&{}+ 24\sigma^2 \log(ep ) \frac{ m \lambda_1(\Sigma) }{\lambda^2} \bigl\{ \EE\bigl[
\widehat\Delta^2_k\bigr] +k\EE\bigl[ d_1^2(P_{\tilde J }E)
\bigr] \bigr\}.
\end{eqnarray*}
Taking $\lambda^2 = C \log(ep) mk \sigma^2 \lambda_1(\Sigma)$ with
$C=792$, we find
\[
64 k \EE\bigl[ d_1^2(P_{\tilde J}E)\bigr] \le
66k \bigl[ 6\sigma^2 n + 16\sigma^2 \exp(-n/2) + 12
\sigma^2 \log(ep ) |J| \bigr] + 2\EE\bigl[\widehat
\Delta^2_k \bigr].
\]
Inserting this bound in (\ref{drie}), we obtain
%
\begin{equation}
\label{DOI} \EE\bigl[ \widehat\Delta_k^2\bigr] \le 10
\Delta^2 +66k \bigl[ 6\sigma^2 n + 16\sigma^2
\exp(-n/2) + 12\sigma^2 \log(ep ) |J| \bigr].
\end{equation}
This concludes the second case. Our risk bound (\ref{twee}) follows
directly from (\ref{UNO}) and~(\ref{DOI}).

\subsection{\texorpdfstring{Proof of Theorem \protect\ref{thmvijf}}{Proof of Theorem 4}}\label{secA.4}
Recall that $\widehat r $ is the number of eigen-values of $Y'PY$ that
exceed the threshold level $2\mu=4\sigma^2 (\sqrt{n}+\sqrt{q})^2$.
Theorem 2 and Corollary 4 in \citet{BuneaSheWegkamp} show that for
$c_0=(\sqrt{2}-1)^2 /2$,
\[
\PP\{ \widehat r \ne r \} \le\exp\bigl\{-c_0(n+q)\bigr\}.
\]
Next, we decompose the risk as follows:
%
\begin{eqnarray}
\label{dc} \EE \bigl[ \bigl\| XA-X\widehat B^{(1)}\bigr\|_F^2
\bigr]
&=& \EE \bigl[ \bigl\| XA-X\widehat B^{(1)}\bigr\|_F^2 I
\{\widehat r=r\} \bigr]\nonumber\\[-8pt]\\[-8pt]
&&{} + \EE \bigl[ \bigl\| XA-X\widehat B^{(1)}
\bigr\|_F^2I\{ \widehat r\ne r\} \bigr].
\nonumber
\end{eqnarray}
The first term on the right gives the bound obtained in Theorem \ref
{thmdrie} for $k=r$. It remains to bound the second term.

Let $O$ denote the $p\times n$ matrix with all entries equal to zero.
Then, since $\widehat r =r(\widehat B) \ge r(O)=0$, we have the inequality
\[
\bigl\| Y-X\widehat B^{(1)}\bigr\|_F^2 + 2\lambda\bigl\|
\widehat B^{(1)}\bigr\|_{2,1} \le \| Y- XO \|_F^2
+ 2\lambda\| O\|_{2,1} = \| Y\|_F^2
\]
by the minimzing property of $\widehat B^{(1)}$.
Using Pythagoras,
$\| PY-X\widehat B^{(1)}\|_F^2 \le\| P Y\|_F^2$ for $P$ the
projection matrix on the column space of $X$.
Consequently,
\begin{eqnarray*}
&& \EE \bigl[ \bigl\| XA-X\widehat B^{(1)}\bigr\|_F^2I\{
\widehat r\ne r\} \bigr]
\\
&&\qquad\le2 \EE \bigl[ \bigl\{ \|P Y- XA\|_F^2 +\bigl\| P Y-X
\widehat B^{(1)}\bigr\|_F^2 \bigr\} I\{ \widehat r\ne
r\} \bigr]
\\
&&\qquad\le2 \EE \bigl[ \bigl\{ \|P E\|_F^2 + \|P Y
\|_F^2 \bigr\} I\{ \widehat r\ne r\} \bigr]
\\
&&\qquad\le6 \EE \bigl[ \| P E \|_F^2 I\{ \widehat r\ne r\}
\bigr] + 4\| XA\|_{F}^2 \PP\{ \widehat r\ne r\}.
\end{eqnarray*}
Since $\| P E\|_F^2/\sigma^2$ has a Chi-square distribution with $qn$
degrees of freedom [see Lemma 3 in \citet{BuneaSheWegkamp}],
we have
\[
\EE \bigl[ \| P E \|_F^4 \bigr] \le
\bigl(q^2 n^2+ 2 q n\bigr) \sigma^4 \le2
q^2 n^2 \sigma^4.
\]
We
obtain, using the Cauchy--Schwarz inequality,
\[
\EE \bigl[ \| P E \|_F^2 I\{ \widehat r\ne r\} \bigr]
\le\sqrt2 qn \sigma^2\exp\bigl\{- c_0 (n+q)/2 \bigr\}
\]
and so
\begin{eqnarray*}
\EE \bigl[ \bigl\| XA-X\widehat B^{(1)}\bigr\|_F^2I\{
\widehat r\ne r\} \bigr]
&\le& 4 \|XA\|_{F}^2 \exp\bigl\{- c_0 (n+q)
\bigr\}\\
&&{} + 6\sqrt2 qn \sigma^2\exp\bigl\{- c_0 (n+q)/2 \bigr\}.
\end{eqnarray*}
The second term on the right is clearly bounded.
For the first term on the right, we invoke condition $\mathfrak{C2}$ on
$\| XA\|_F^2$. This completes our proof.

\subsection{\texorpdfstring{Proof of Theorem \protect\ref{thmvier}}{Proof of Theorem 5}}\label{secA.5}
In this proof, $C_1,\ldots, C_7$ are numerical, positive and finite constants.
Theorem~\ref{thmtwee} of Section~\ref{sec2}, applied to the random matrices
$\widehat B_{\lambda,k}$, yields
\begin{eqnarray*}
&& \EE\bigl[ \bigl\| X\widehat B^{(2)}- XA\bigr\|_F^2
\bigr]
\\
&&\qquad \le\min_{k, \lambda} \EE \bigl[ 10 \|X\widehat B_{\lambda,k} -XA
\|_F^2 + 8 \operatorname{pen} (\widehat
B_{\lambda,k}) \bigr] + 768 n\sigma^2 \exp(-n/2)
\\
&&\qquad \le C_1 \min_\lambda\EE \bigl[ \|X\widehat
B_{\lambda,r} -XA\|_F^2 + \operatorname{pen} (
\widehat B_{\lambda,r}) +n \sigma^2 \bigr].
\end{eqnarray*}
(Here we used $c=12$ in the penalty term.)
Since $\Sigma$ satisfies $\mathfrak{A}(J(A),\delta_{J(A)})$, and we
assume that $\lambda_1(\Sigma)/ \delta_{J(A)} \le C_3$,
Theorem~\ref{thmdrie} yields, for each global solution~$\widehat
B_{\lambda,r}$,
\[
\EE\bigl[\| X\widehat B_{\lambda,r} -XA\|_F^2
\bigr] \le C_2 \bigl\{ n+\bigl|J(A)\bigr|\log(p) \bigr\} r \sigma^2
\]
for $\lambda^2= C\sigma^2 \lambda_1(\Sigma) k m \log(ep)$ with $C$
large enough.
It remains to bound the expected penalty term $\EE[\operatorname
{pen}(\widehat B_{\lambda,r})]$. Since
\[
\EE\bigl[\operatorname{pen}(\widehat B_{\lambda,r})\bigr] \le
C_4 \bigl( n+ \log (p)\EE\bigl[ \bigl| J(\widehat B_{\lambda,r}) \bigr|
\bigr] \bigr),
\]
we need to bound $\EE[ | J(\widehat B_{\lambda,r}) | ]$. We write
$\widehat J_r= J(\widehat B_{\lambda,r})$.
From (\ref{vijf}) in the proof of Theorem~\ref{thmdrie}, we have
\[
\EE\bigl[ |\widehat J_r | \bigr] \le \frac{2m \lambda_1(\Sigma)}{\lambda^2} \bigl\{ \EE\bigl[
\| X\widehat B_{\lambda,r} -XA\|_F^2 \bigr] + r \EE
\bigl[ d_1^2 (P_{\widehat J_r}E)\bigr] \bigr\},
\]
while Lemma~\ref{maxbound} gives
\[
\EE\bigl[ d_1^2(P_{\widehat J_r } E )\bigr] \le 6n
\sigma^2 + 16\sigma^2 \exp (-n/2) + 12
\sigma^2 \log(ep) \EE\bigl[ |\widehat J|\bigr].
\]
Therefore, taking $\lambda^2= C\sigma^2 \lambda_1(\Sigma) k m \log(ep)$
with $C$ large enough, we obtain from the previous three displays
\begin{eqnarray*}
\EE\bigl[ |\widehat J_r| \bigr] &\le& C_5 \bigl\{ \EE\bigl[\| X
\widehat B_{\lambda,r} -XA \|_F^2 \bigr] / r + n
\sigma^2 / \log(p) \bigr\}
\\
&\le& C_6 \sigma^2 \bigl\{ n+\bigl|J(A)\bigr|\log(p) \bigr\}.
\end{eqnarray*}
We now can conclude that
\[
\EE\bigl[ \bigl\| X\widehat B^{(2)}- XA\bigr\|_F^2
\bigr] \le C_7 \sigma^2 r \bigl\{ n+ \log(p) \bigl|J(A)\bigr|
\bigr\},
\]
and the proof is complete.

\subsection{\texorpdfstring{Proof of Theorem \protect\ref{algstat}}{Proof of Theorem 6}}\label{secA.6} \label
{subsecglobalconv}

In this proof we drop the subscripts $\lambda, k$ in the $(S, V)$
iterates.\vspace*{6pt}

\textsc{Proof of part} (ii). The proof of this part of Theorem \ref
{algstat} will follow
from the global convergence theorem (GCT) of \citet{zangwill69}. For
completeness, we state
this theorem below, then we verify that its conditions hold in Lemmas
\ref{algconv},~\ref{unifbdd} and~\ref{closedmap}.\vspace*{-2pt}

\begin{theo*}[{[\citet{Luenberger08}]}]
Let $A$ be a map describing an algorithm on $\mathcal{X}$ and suppose
that given $x_0$ the sequence $\{x_k\}_{k=1}^{\infty}$ is generated by
$x_{k+1} \in A(x_k)$.
Let a solution set $\Gamma\subset\mathcal{X}$ be given, and suppose:
\begin{longlist}[(2)]
\item[(1)] All points $x_k$ are contained in a compact set $S\subset
\mathcal{X}$;
\item[(2)] There exists a continuous function $Z$ on $\mathcal{X}$ such
that \textup{(a)} if $x\notin\Gamma$, then $Z(y)< Z(x)$ for all $y \in
A(x)$;
\textup{(b)} if $x\in\Gamma$, then $Z(y) \leq Z(x)$ for all $y\in A(x)$;
\item[(3)] The mapping $A$ is closed at points outside $\Gamma$.
\end{longlist}
Then the limit of any convergent subsequence of $\{x_k\}_k$ is a
solution.\vspace*{-2pt}
\end{theo*}

We begin by introducing a map, usually referred to in the literature as
a \textit{point-to-set} map, to characterize our Algorithm \ref
{algS-RRR}. Let $\Omega= \mathbb{R}^{p\times k} \times\mathbb
{O}^{n\times k}$.
Define ${\mathcal M}^S\dvtx  \Omega\rightarrow2^\Omega, {\mathcal M}^V\dvtx
\Omega\rightarrow2^\Omega$ as follows:
\begin{eqnarray*}
{\mathcal M}^S(S, V) &=& \Bigl\{(\bar S, V)\in\Omega\dvtx
\inf_{\tilde S\in
\mathbb{R}^{p\times k}} F(\tilde S, V) \geq F(\bar S, V)\Bigr\},
\\[-2pt]
{\mathcal M}^V(S, V) &=& \Bigl\{(S, \bar V)\in\Omega\dvtx
\inf_{\tilde V \in
\mathbb{O}^{n\times k}} F(S, \tilde V) \geq F(S, \bar V)\Bigr\}
\end{eqnarray*}
and define $\mathcal{M}={\mathcal M}^V {\mathcal M}^S$ as a \textit
{composite} point-to-set map; see \citet{Luenberger08} for more details.
Algorithm~\ref{algS-RRR} can be described by
\[
\bigl(S^{(j+1)}, V^{(j+1)}\bigr)\in\mathcal{M}
\bigl(S^{(j)}, V^{(j)}\bigr),
\]
that is, $(S^{(j+1)}, V^{(j)})\in{\mathcal M}^S(S^{(j)}, V^{(j)})$ and
\mbox{$(S^{(j+1)}, V^{(j+1)})\in{\mathcal M}^V(S^{(j+1)}, V^{(j)})$}.

Recall that
\[
F(S, V; \lambda) =\tfrac{1}{2} \bigl\| Y - X S V'
\bigr\|_F^2 + \lambda\bigl\| S V' \bigr\|_{2,1},\qquad
S\in\mathbb{R}^{p\times k}, V\in\mathbb{O}^{n\times k},
\]
and that we analyze the \textit{unconstrained} minimum of $F$ over the
product manifold $\mathbb{R}^{p\times k} \times\mathbb{O}^{n\times
k}$. For simplicity, we will write just $F(B)$ and $F(S, V)$ for
$F(B;\lambda)$ and $F(S, V;\lambda)$, respectively, when there is no
ambiguity.

Lemma~\ref{algconv} shows the algorithm converges globally for any
initial starting point.\vspace*{-2pt}
%
\begin{lemma} \label{algconv}
For any $j\geq0$, $\operatorname{rank}(B^{(j)}) \leq k$ and
$F(B^{(j)}) \geq F(B^{(j+1)})$.\vspace*{-2pt}
\end{lemma}
\begin{pf}
We write $F(S, V)=\frac{1}{2} \| Y V - XS\|_F^2 + \lambda\| S\|_{2,1}+\frac{1}{2}
\operatorname{tr}(Y (I -\break V V') Y')$.
Given $V$, (\ref{grpsel-constrRank2}) reduces to\vadjust{\goodbreak} the following
optimization problem after vectorization:
%
\begin{equation}\label{subprob1}
\tfrac{1}{2} \bigl\|{\operatorname{vec}}\bigl( (Y V)'\bigr) - ( X \otimes I)
\operatorname{vec}\bigl( S'\bigr)\bigr\|_F^2 + \lambda
\|S \|_{2,1},
\end{equation}
where ``$\mathrm{vec}$'' is the standard vectorization operator
and $\otimes$ is the Kronecker product. This is a GLASSO-type
optimization problem that is convex in $S$. The global minimum of (\ref
{subprob1}) can always be achieved at some $S$ with $\| S \|<\infty$.
Given $S$, writing $F(S, V)=- \operatorname{tr}(Y' X S V') + \|Y\|_F^2/2 + \| X S \|_F^2/2 + \lambda\| S \|_{2,1}$, we see the
optimization problem is equivalent to
%
\begin{equation}\label{subprob2}
\max_{V\in{\mathbb O}^{n\times p}} \operatorname{tr}\bigl(W' V\bigr),
\end{equation}
where $W = Y' X S \in\mathbb{R}^{n\times k}$. The (global) maximum can
be attained, too, due to the compactness of ${\mathbb O}^{n\times p}$.
In fact, $ \operatorname{tr}(W' V) \leq\sum_i d_i(W)$ by von Neumann's
trace inequality.
Let $W=U_w D_w V_w'$ be the SVD with $V_w\in\mathbb{O}^{k\times k}$. Then
%
\begin{equation}\label{optV}
\widehat V = U_w V_w'
\end{equation}
achieves the upper bound $\sum_i d_i(W)$. This globally optimal
solution to (\ref{subprob2}) is the one used in Algorithm~\ref{algS-RRR}.
Therefore, we have
\[
F\bigl(V^{(j)}, S^{(j)}\bigr) \geq F\bigl(V^{(j)},
S^{(j+1)}\bigr) \geq F\bigl(V^{(j+1)}, S^{(j+1)}\bigr),
\]
that is, $F(B_{k,\lambda}^{(j)}) \geq F(B_{k,\lambda}^{(j+1)})$ during
each iteration.\vspace*{-2pt}
\end{pf}
%
\begin{lemma} \label{unifbdd}
Suppose $\lambda>0$. Then $B^{(j)}$, $S^{(j)}$, $V^{(j)}$ in
Algorithm~\ref{algS-RRR} are uniformly bounded in $j$.\vspace*{-2pt}
\end{lemma}
\begin{pf}
From Lemma~\ref{algconv}, $ \|S^{(j)}\|_{2,1} \leq
F(S^{(j)}, V^{(j)})/\lambda\leq F(S^{(1)}, V^{(1)})/\lambda\leq
F(0, V^{(0)})/\lambda\le\| Y \|_F^2/(2\lambda)$. Therefore, $\|
S^{(j)}\|$ must be uniformly bounded.\vspace*{-2pt}~%
\end{pf}
%
\begin{lemma} \label{closedmap}
Suppose $\lambda>0$. The map $\mathcal{M}$ introduced for describing
Algorithm~\ref{algS-RRR} is a closed point-to-set map on $\Omega$.\vspace*{-2pt}
\end{lemma}
\begin{pf}
Notice that:
\begin{itemize}
\item The set $\Omega=\mathbb{R}^{p\times k} \times\mathbb
{O}^{n\times
k}$ is closed because $\mathbb{O}^{n\times k}$ is the inverse image of
$\{I\}$ under the continuous function $ M \mapsto M^T M$, $\forall M
\in\mathbb{R}^{n\times k}$; $\mathbb{O}^{n\times k}$ is in fact an
embedded submanifold of $\mathbb{R}^{n\times k}$.
\item${\mathcal M}^S(\omega)\neq\varnothing, {\mathcal M}^V(\omega
)\neq\varnothing, \forall\omega\in\Omega$, seen from the proof of
Lemma~\ref{algconv}.
\end{itemize}
First, we prove that ${\mathcal M}^S$ is closed on $\Omega$.
It suffices to show the point-to-set map ${\mathcal M}^s (x) =\{y \in
{\mathbb R}^{p\times k}\dvtx  F(y, x) \leq\min_{\tilde y\in{\mathbb
R}^{p\times k}} F(\tilde y, x) \}$ is closed at any $x \in{\mathbb
O}^{n\times k}$.
Let $x_j \rightarrow x^*$, $y_j \in{\mathcal M}^s(x_j)$ and
$y_j\rightarrow y^*$, with $x_j, x^* \in\mathbb{O}^{n\times k}, y_j,
y^*\in\mathbb{R}^{p\times k}$. Suppose $y^* \notin{\mathcal
M}^s(x^*)$. Since $y^* \in{\mathbb R}^{p\times k}$, $F(y^*, x^*) >
\min_{y\in{\mathbb R}^{p\times k}} F(y,\break x^*) =: L$. There exists some
$\varepsilon_0>0$ such that $F(y^*, x^*) > L+\varepsilon_0$.
Let $\tilde y \in{\mathcal M}^s(x^*)$. Then $F(\tilde y, x^*)=L$.
Since $F$ is continuous at $(\tilde y, x^*)$, $\lim_{j\rightarrow
\infty
} F(\tilde y,\break x_j)=L$. For $j$ large enough, $| F(\tilde y, x_j) - L|
\leq\varepsilon_0/2$, from which it follows that $F(y_j, x_j) \leq
F(\tilde y, x_j) \leq L+\varepsilon_0/2$ and $F(y^*, x^*) \leq
L+\varepsilon_0/2$.\vadjust{\goodbreak}

The contradiction implies $y^* \in{\mathcal M}^s(x^*)$.
Hence, ${\mathcal M}^{s}$ and thus ${\mathcal M}^{S}$ are closed.
From the proof of Lemma~\ref{unifbdd}, we also know ${\mathcal M}^S$
is compact.
Similarly, we can show ${\mathcal M}^V$ is closed on ${\mathbb
R}^{p\times k}$. Based on the properties of the point-to-set maps
[\citet{Luenberger08}, page 205], $\mathcal{M}$ is closed on
$\Omega$.
\end{pf}

If we set $Z$ in the general statement of the global convergence
theorem to be our continuous criterion function $F$, and if we take
$\Gamma$ to be the set of local of minima of~$F$, Lemma~\ref{algconv}
and the assumption in part (ii) of our Theorem~\ref{algstat2} guarantee that (2) of
GCT holds.
Lemmas~\ref{unifbdd} and~\ref{closedmap} verify conditions (1) and (3)
of GCT, respectively.
This concludes the proof of this part.\vspace*{8pt}

\textsc{Proof of part} (i).
From displays (\ref{subprob1}) and (\ref{subprob2}) we observe that $F$
is \textit{convex} in $S$, given $V$, and it is linear and therefore
\textit{smooth} in $V$, given $S$. This part of the theorem shows that,
under no further conditions on
$F$, Algorithm~\ref{algS-RRR} converges to a stationary point of $F$.
We begin by defining
a stationary point in this context. Recall that we view (\ref
{grpsel-constrRank2}) as an unconstrained optimization problem in
$\Omega
={\mathbb R}^{p\times k} \times{\mathbb O}^{n\times k}$.
Notice that $\mathbb{O}^{n\times k}$, which is a Stiefel manifold,
is a Riemannian submanifold of $\mathbb{R}^{n\times k}$. We use the
inherited Riemannian metric to define the gradient of $F$ with respect
to $V$ [\citet{Boothby02}]. This Riemannian gradient, denoted by
$\nabla_V F(S, V)$, can be explicitly computed: $\nabla_V F(S, V) = {\mathcal
P}_V (-W) = -W + V V' W/2 + V W' V/2$ with $W=Y' X S$, where ${\mathcal
P}_V$ is the projection onto the tangent space to ${\mathbb O}^{n\times
k}$ at $V$.

Since $F$ is convex in $S$, the subdifferential of $F$ with respect to
$S$, denoted by $\partial_S F$, is also well defined. From \citet
{Gabay82} and Shimizu, Ishizuka and Bard [(\citeyear{Shim}), page 62],
a necessary condition for $F$ to have a local minimum at $(S^*, V^*)$
is that this point is stationary, that is,
%
\begin{equation}\label{graddef}
0 \in\partial_S F\bigl(S^*, V^*\bigr)\quad \mbox{and}\quad
\nabla_V F\bigl(S^*, V^*\bigr) = 0.
\end{equation}

For future use, note that $\nabla_V F(S, V)$ is continuous on $\Omega$.
Since $V^{(j + 1)}$ minimizes $F(S^{(j+1)}, \cdot)$, we have $\nabla_V
F(S^{(j+1)}, V^{(j+1)})=0$, for any $j\geq0$.

Because the optimum of (\ref{subprob1}) [or (\ref{subprob2})] may not
be uniquely attained, $F$ is not guaranteed to be a \textit{descent
function} in general [see \citet{zangwill69}, \citet{Bert}
and \citet
{Luenberger08} for details]. Therefore, GCD cannot be directly applied.

From Lemma~\ref{algconv}, $F(S^{(j)}, V^{(j)})$ must converge. Denote
the limit by $L^*$. Let $(S^{(j_l)}, V^{(j_l)})$ ($l=1, 2, \ldots$) be
a subsequence of $(S^{(j)}, V^{(j)})$ which converges to $(S^{*},
V^{*})$ as $l\rightarrow\infty$. Then $L^*=\lim_{l\rightarrow\infty}
F(S^{(j_l)}, V^{(j_l)}) = F(S^{*}, V^{*})$. We assume, without loss of
generality, $(S^{(j_l+1)}, V^{(j_l+1)})$ also converges (Lemma~\ref
{unifbdd}), and denote the limit by $(\bar S, \bar V)$. We have
%
\begin{equation}\label{contr0}
F(\bar S, \bar V)= L^*.
\end{equation}

We claim that $(S^*, V^*)$ must be a stationary point of $F$.
First, the continuity of
$\nabla_V F(S, V)$ and the fact that $\nabla_V F(S^{(j)}, V^{(j)})=0$
imply that \mbox{$\nabla_V F(S^*, V^*)=0$}.\vadjust{\goodbreak} Suppose $0 \notin\partial_S
F(S^*, V^*)$. This implies $S^*$ is not a global minimizer of (\ref{subprob1}).
Lemma~\ref{closedmap}, however, states that $(\bar S, \bar V) \in
{\mathcal M}(S^{*}, V^{*})$.
That is, $ L^* = F(\bar S, \bar V) \leq F(\bar S, V^*) < F(S^*, V^*) =
L^*$, which contradicts (\ref{contr0}). The last inequality is strict
because $F$ is convex in $S$ and so, applying the algorithm to $S^*$,
which is not a global minimizer, yields a strict improvement (decrease)
of the criterion function. Hence, $0 \in\partial_S F(S^*, V^*)$. The
proof is complete.\vspace*{-2pt}

\subsection{\texorpdfstring{Proof of Theorem \protect\ref{algstat2}}{Proof of Theorem 7}}\label{secA.7}
\label{subsecglobalconv2}

We use the same notation system as in Appendix~\ref{subsecglobalconv}.
Recall that $(S^*, V^*)\in\Omega$ is a coordinatewise minimum point of
$F(S, V)$ if $F(S, V^*)\geq F(S^*, V^*)$, $\forall S\in{\mathbb
R}^{p\times k}$ and $F(S^*, V)\geq F(S^*, V^*)$, $\forall V\in{\mathbb
O}^{n\times k}$ [\citet{Tseng01}]. In our problem, this implies $(S^*,
V^*)$ is also a stationary point of $F(S, V)$.

Without loss of generality, assume $X$ has been scaled to have $\|X\|_2\leq1$ before running algorithm $\mathscr A'$.
For simplicity, set $K=1$ in (\ref{grpsel-th}) and redefine the
operator $T_V\dvtx  {\mathbb R}^{p\times k} \rightarrow{\mathbb R}^{p\times
k}$ by
%
\begin{equation}\label{Tdef2}
T_V \circ S = \vec\Theta \bigl( X' Y V + \bigl(I -
X' X\bigr) S; \lambda \bigr)\qquad \forall S \in{\mathbb
R}^{p\times k}.
\end{equation}
Let $T_V^\alpha$ ($\alpha\in\mathbb N$) be the composition of
$\alpha$ $T_V$'s.
Define point-to-set maps ${\mathcal M}^S(S, V) = \{(\bar S,\bar V)\in
\Omega\dvtx  \bar S \in\{ T_V \circ S, T_V^2 \circ S,\ldots,
T_V^{M_{\mathrm{iter}}} \circ S\}, \bar V=V\}$,
and ${\mathcal M}^V(S, V) = \{(\bar S, \bar V)\in\Omega\dvtx  \inf_{\tilde
V \in\mathbb{O}^{n\times k}} F(S, \tilde V) \geq F(S, \bar V), \bar
S=S\}$.
Then $\mathcal{M}={\mathcal M}^V {\mathcal M}^S$ characterizes
$\mathscr A'$. When updating $S$ at step (a), the algorithm allows one
to perform $T$ any times (denoted by $\alpha_j$) provided $\alpha_j$
does not go beyond $M_{\mathrm{iter}} \in\mathbb N$ that is prespecified before
running $\mathscr A'$.\vspace*{-2pt}
%
\begin{lemma} \label{senstS}
Given any $V\in{\mathbb O}^{n\times k}$ and $S\in{\mathbb R}^{p\times
k}$, let $\tilde S = T_V(S)$. Then
$F(S, V) - F(\tilde S, V) \geq\frac{1}{2} \| \tilde S - S \|_F^2$.\vspace*{-2pt}
\end{lemma}
\begin{pf}
Apply the theorem in \citet{SheGLMTISP} to the vectorized problem
(\ref
{subprob1}). Note that $2 - \| X \otimes I \|_2 = 2- \| X \|_2 \geq1$.
The proof details are omitted.\vspace*{-2pt}~%
\end{pf}
Choose $(\bar S, V) \in{\mathcal M}^S(S, V)$, using the triangle
inequality we know
$F(S, V) - F(\bar S, V) \geq\frac{1}{2} \| \bar S - S \|_F^2$.\vspace*{-2pt}
%
\begin{lemma} \label{unifbdd2}
Suppose $\lambda>0$. Then $B^{(j)}$, $S^{(j)}$, $V^{(j)}$ in $\mathscr
A'$ are uniformly bounded in $j$.\vspace*{-2pt}
\end{lemma}
The proof is similar to the proof of Lemma~\ref{unifbdd} and therefore
omitted.\vspace*{-2pt}
%
\begin{lemma}
Suppose $\lambda> 0$.
The $\mathcal{M}$ for describing $\mathscr A'$ is a closed point-to-set
map on $\Omega$.\vspace*{-2pt}
\end{lemma}
\begin{pf}
Similar to the proof of Lemma~\ref{closedmap}, we prove that the
point-to-set map ${\mathcal M}^s (S, V) =\{ T_V \circ S, T_V^2 \circ S,\ldots,
T_V^{M_{\mathrm{iter}}} \circ S\}$ is closed at any $(S, V) \in\Omega$.
Then $M^S$ and thus ${\mathcal M}$ are closed on $\Omega$.

Let $(S_j, V_j) \rightarrow(S^*, V^*)$, $\tilde S_j \in{\mathcal
M}^s(S_j, V_j)$ and $\tilde S_j\rightarrow\tilde S^*$, with $(S_j,
V_j)$, $(S^*,\break V^*) \in\Omega$ and\vadjust{\goodbreak} $\tilde S_j, \tilde S^*\in
\mathbb{R}^{p\times k}$.
There must exist infinitely many $\tilde S_j$'s satisfying
$T_{V_j}^{\alpha} \circ(S_j, V_j) = \tilde S_j$ for some $\alpha\in
{\mathbb N}$. Let $g(S, V) = T_V^{\alpha}(S, V)$. It is not difficult
to see that $g$ is jointly continuous. Hence, $\tilde S^* \in{\mathcal
M}^s(S^*, V^*)$ by a subsequence argument.\vspace*{-2pt}
\end{pf}

Now we prove Theorem~\ref{algstat2}.
Following the lines of the proof of part (i) of Theorem~\ref{algstat},
for any accumulation point $(S^*, V^*)$ of $(S^{(j)}, V^{(j)})$, $(S^*,
V^*)\in\Omega$ and there exists $(\bar S, \bar V)\in{\mathcal M}(S^*,
V^*)$ with $F(\bar S, \bar V) = F(S^*, V^*)$. Since $F(\bar S, \bar
V)\leq F(\bar S, V^*) \leq F(S^*, V^*)$, $F(\bar S, V^*) = F(S^*,
V^*)$. It follows from the comment after Lemma~\ref{senstS} that $\bar
S = S^*$.
This means $T_{V^*}^{\alpha_0} \circ S^* = S^*$ for some $\alpha_0\in
{\mathbb N}$. But then $F(T_{V^*}^{\alpha} \circ S^*, V^*) = F(S^*,
V^*)$ for any $\alpha\leq\alpha_0$, and, in particular, $F(T_{V^*}
\circ S^*, V^*) = F(S^*, V^*)$.
Applying Lemma~\ref{senstS} again yields
$T_{V^*} \circ S^* = S^*$. It is easy to verify from (\ref{Tdef2}) that
$S^*$ is a fixed point of $T_{V^*}$ is equivalent to $0\in\partial_S
F(S^*, V^*)$. Therefore, $S^*$ is a global minimizer of $F(S, V^*)$
given $V^*$, due to the convexity of~(\ref{subprob1}).

On the other hand, from $\bar S = S^*$, we have $(S^*, \bar V) \in
{\mathcal M}^V(S^*, V^*)$. $\bar V$ is a (global) minimizer of $F(S^*,
V)$ given $S^*$. But $F(S^*, \bar V) = F(\bar S, \bar V) = F(S^*,
V^*)$, so $V^*$ also minimizes $F(S^*, V)$ given $S^*$, and $\nabla_V
F(S^*, V^*)=0$.
In summary, $(S^*, V^*)$ is a coordinatewise minimum of $F$.\vspace*{-2pt}
\end{appendix}


%

\printaddresses

\end{document}